\documentclass[a4paper,11pt,twoside]{article}
\usepackage{mfpaperstuff2}
\usepackage{ifthen}

\makeatletter
\newsavebox{\@brx}
\newcommand{\llangle}[1][]{\savebox{\@brx}{\(\m@th{#1\langle}\)}%
  \mathopen{\copy\@brx\kern-0.7\wd\@brx\usebox{\@brx}}}
\newcommand{\rrangle}[1][]{\savebox{\@brx}{\(\m@th{#1\rangle}\)}%
  \mathclose{\copy\@brx\kern-0.7\wd\@brx\usebox{\@brx}}}
\makeatother

\newgeometry{margin=0.421in}
\begin{document}
\hyphenation{Rouq-uier}
\newcommand\hkl[3]{h_{#1,#2}(#3)}
\newlength\oneht\settoheight\oneht{$1$}
\newcommand\tim{\begin{tikzpicture}[scale=1,baseline=0pt]\draw(0,.5*\oneht)node{$\times$};\end{tikzpicture}}
\newcommand\bk[2]{B_{#2}}
\newcommand\bks[1]{\bk\sigma{#1}}
\renewcommand\dim[1]{\operatorname{deg}\spn{#1}}
\newcommand\gry{\Yfillcolour{gray}}
\newcommand\wht{\Yfillcolour{white}}
\newcommand\ol[1]{\widebar{#1}}
\newcommand\sse[1]{[#1]}
\newcommand\indi[1]{\mathbbm1(#1)}
\newcommand\lad[2]{\operatorname{lad}_{#1}(#2)}
\newcommand\alad[2]{\operatorname{lad}^+_{#1}(#2)}
\newcommand\rlad[2]{\operatorname{lad}^-_{#1}(#2)}
\newcommand\slp[2]{\operatorname{slp}_{#1}(#2)}
\newcommand\aslp[2]{\operatorname{slp}^+_{#1}(#2)}
\newcommand\rslp[2]{\operatorname{slp}^-_{#1}(#2)}
\newcommand\zslp[2]{\operatorname{slp}^\ast_{#1}(#2)}
\newcommand\noregdown[2]{#2\lightning_{#1}}
\newcommand\ddim[1]{\operatorname{ddeg}\spn{#1}}
\newcommand\df{\stackrel4\dots}
\newcounter{case}\setcounter{case}0
\newcommand\nextcase{\addtocounter{case}1Case \arabic{case}: }
\newcommand\dn[1]{{\downarrow}_{#1}}
\newcommand\up[1]{{\uparrow}^{#1}}
\renewcommand\labelitemii{$\diamond$}
\newcommand\fs[3]{G(#1,#2,#3)}
\newcommand\strp[2]{\operatorname{Hs}_{#1\setminus#2}}
\newcommand\lrcb[2]{\lrc{#1}{#2}{(1^\bullet)}}
\newcommand\quo[2]{#1^{(#2)}}
\newcommand\epstar[1]{\epsilon^\ast(#1)}
\newcommand\ddd{\partial}
\newcommand\spe[1]{\mathrm S^{#1}}
\newcommand\jms[1]{\mathrm D^{#1}}
\newcommand\schir[1]{\mathrm L^{#1}}
\newcommand\weyl[1]{\Delta^{#1}}
\newcommand\sym{\mathrm{Sym}}
\newcommand\cm[2]{[#1:#2]}
\newcommand\muset[2]{\{0^{#1},1^{#2}\}}
\newcommand\ev[1]{\operatorname{ev}(#1)}
\newcommand\spre{spin-removable\xspace}
\newcommand\sprm{spin-removable node\xspace}
\newcommand\sprms{spin-removable nodes\xspace}
\newcommand\spam{spin-addable node\xspace}
\newcommand\spams{spin-addable nodes\xspace}
\newcommand\ird\scri
\newcommand\treg{$2$-regular\xspace}
\newcommand\trp{\treg partition\xspace}
\newcommand\trps{\treg partitions\xspace}
\newcommand\fbc{$4$-bar-core\xspace}
\newcommand\fbcs{$4$-bar-cores\xspace}
\newcommand\bc{bar-core\xspace}
\newcommand\schur[1]{\cals(#1)}
\newcommand\qschur[1]{\cals^\circ(#1)}
\newcommand\odc{\mathring D}
\newcommand\dspn{D^{\mathsf{spn}}}
\newcommand\decs[2]{D^{\mathsf{spn}}_{#1#2}}
\renewcommand\v{^{-1}}
\newcommand\ip[2]{\lan#1,#2\ran}
\newcommand\qs{separated\xspace}
\newcommand\Qs{Separated\xspace}
\newcommand\rs[2]{#1{\downarrow}_{#2}}
\newcommand\spr{spin residue\xspace}
\newcommand\sprs{spin residues\xspace}
\newcommand\dup[1]{#1\sqcup#1}
\newcommand\chm[2]{\left(#1\!:\!#2\right)}
\renewcommand\deg{\operatorname{deg}}
\renewcommand{\crefrangeconjunction}{--}
\medmuskip=3mu plus 1mu minus 1mu
\renewcommand{\rt}[1]{\rotatebox{90}{$#1$}}
\newcommand\aaa{\mathfrak{A}_}
\newcommand\taaa{\tilde{\mathfrak{A}}_}
\newcommand\tsss{\tilde{\mathfrak{S}}_}
\newcommand\len[1]{l(#1)}
\newcommand\res{{\downarrow}}
\newcommand\ind{{\uparrow}}
\newcommand\reg{^{\operatorname{reg}}}
\newcommand\dbl{^{\operatorname{dbl}}}
\newcommand\dblreg{^{\operatorname{dblreg}}}
\newcommand\ee[1]{\mathrm{e}_{#1}}
\newcommand\ff[1]{\mathrm{f}_{#1}}
\newcommand\eed[2]{\ee{#1}^{(#2)}}
\newcommand\ffd[2]{\ff{#1}^{(#2)}}
\newcommand\emx[1]{\eed{#1}{\operatorname{max}}}
\newcommand\fmx[1]{\ffd{#1}{\operatorname{max}}}
\newcommand\ord[1]{\llbracket#1\rrbracket}
\newcommand\spn[1]{\llangle#1\rrangle}
\newcommand\spq[3]{\spn{#1+4#2\sqcup2#3}}
\newcommand\sid[1]{\varphi(#1)}
\newcommand\prj[1]{\operatorname{prj}(#1)}
\newcommand\modr[1]{\widebar{#1}}
\newcommand\mspn[1]{\modr{\spn{#1}}}
\newcommand\ads[2]{\ifthenelse{\equal{#2}1}{\stackrel{#1}\Longrightarrow}{\stackrel{#1^#2}\Longrightarrow}}
\newcommand\lrc[3]{\operatorname{a}^{#1}_{#2#3}}
\newcommand\sstd[1]{\operatorname{Sstd}_X(#1)}
\newcommand\shtd[1]{\operatorname{Shtd}_X(#1)}
\newcommand\shyng[1]{\operatorname{shYng}(#1)}

\Yvcentermath0
\Yboxdim{12pt}

\pdfpagewidth=12.5cm
\pdfpageheight=8.839cm
\setcounter{page}0
{\footnotesize This is the author's version of a work that was accepted for publication in \textit{Advances in Mathematics}. Changes resulting from the publishing process, such as peer review, editing, corrections, structural formatting, and other quality control mechanisms may not be reflected in this document. Changes may have been made to this work since it was submitted for publication. A definitive version was subsequently published in\\
\textit{Adv.\ Math.} \textup{\textbf{374} (2020) 107340.\\
http://dx.doi.org/10.1016/j.aim.2020.107340}\normalsize}
\newgeometry{margin=1in,includehead,includefoot}

\title{Irreducible projective representations of the alternating group which remain irreducible in characteristic $2$}
\runninghead{Irreducible projective representations of the alternating group in characteristic $2$}

\toptitle

\begin{abstract}
For any finite group $G$ it is an interesting question to ask which ordinary irreducible representations of $G$ remain irreducible in a given characteristic $p$. We answer this question for $p=2$ when $G$ is the proper double cover of the alternating group. As a key ingredient in the proof, we prove a formula for the decomposition numbers in Rouquier blocks of double covers of symmetric groups, in terms of Schur P-functions.
\end{abstract}

\pdfpagewidth=210mm
\pdfpageheight=297mm

\tableofcontents

\renewcommand\baselinestretch{1.089}

\section{Introduction}

An important problem in the modular representation theory of finite groups is to determine, for a given group $G$ and prime $p$, which ordinary irreducible representations of $G$ remain irreducible in characteristic $p$. This problem was solved for the symmetric groups in a series of papers \cite{jm1,jmp2,slred,mfred,mfirred}. More recently, the author \cite{mfonaltred,mfaltred} completed the same task for the alternating groups, and in \cite{mfspin2} for the double covers $\tsss n$ of the symmetric groups in characteristic $2$. The main result of the present paper (\cref{mainthm}) is a solution to the same problem for the double covers of the alternating groups in characteristic $2$. We hope to address the odd-characteristic case in future work.

As a more general problem, one can ask whether a given ordinary irreducible representation is \emph{homogeneous} in characteristic $p$; that is, all the composition factors of its reduction modulo $p$ are isomorphic. In fact, it is shown in \cite{mfspin2} that solving our main problem for the double cover $\taaa n$ of the alternating group is equivalent to determining which irreducible representations of $\tsss n$ labelled by partitions with exactly two non-zero even parts are homogeneous in characteristic $2$. So we work almost entirely with $\tsss n$ in this paper, translating our results to $\taaa n$ as the final step. The problem of determining all irreducible representations of $\tsss n$ which are homogeneous in characteristic $2$ remains unsolved, and we hope to address this in future work.

In the course of our work, we prove another substantial result on the decomposition numbers for $\tsss n$. A certain class of blocks of symmetric groups (known as \emph{Rouquier blocks} or \emph{RoCK blocks}) are particularly well understood, and have been an important tool in proving various results relating to the symmetric groups. In characteristic $2$, blocks of $\sss n$ naturally correspond to blocks of $\tsss n$, so one can extend the definition of Rouquier blocks to $\tsss n$. The proof of our main result depends on finding an explicit formula (\cref{mainrouqdec}) for the decomposition numbers of Rouquier blocks of $\tsss n$.

Proving our formula for decomposition numbers in Rouquier blocks requires some work with symmetric functions -- specifically, Schur P-functions, which have long been associated with projective representations of symmetric groups. We take a detour into the theory of symmetric functions to prove an auxiliary formula (\cref{mainsym}) needed for our main result on Rouquier blocks.

We now summarise the layout of this paper. \cref{partnsec} contains a brief review of the combinatorics of partitions needed for our work on both symmetric functions and representation theory. In \cref{symfnsec} we prove the results we need on symmetric functions, briefly introducing background as needed. In \cref{mainthmsec} we introduce the background we need on representations of $\sss n$, $\tsss n$ and $\taaa n$, and state our main theorem. In \cref{rouqsec} we introduce Rouquier blocks, state and prove a formula for their decomposition numbers, and derive the consequences for our main theorem. In \cref{mainproofsec} we complete the proof of our main theorem by induction. We end with an index of notation.

\begin{acks}
The research in this paper would not have been possible without extensive computations using GAP \cite{gap}. The author would also like to express his gratitude to an anonymous referee for a very thorough reading of the paper and numerous excellent suggestions for improvements to the text.
\end{acks}

\section{Background on partitions}\label{partnsec}

\subsection{Compositions and partitions}\label{compnpartnsec}

A \emph{composition} is defined to be an infinite sequence $\la=(\la_1,\la_2,\dots)$ of non-negative integers with finite sum. If $\la$ is a composition, its \emph{size} is the sum $|\la|=\la_1+\la_2+\dots$, and we say that $\la$ is a composition of $|\la|$. The integers $\la_1,\la_2,\dots$ are referred to as the \emph{parts} of $\la$.

A \emph{partition} is a composition $\la$ such that $\la_1\gs\la_2\gs\dots$. We write $\calp$ for the set of all partitions. If $\la\in\calp$, we say that $\la$ is a partition of $|\la|$. We write $\calp(n)$ for the set of all partitions of $n$. We define $\len\la$ to be the largest $r$ for which $\la_r>0$, and we say that $\la$ has \emph{length} $\len\la$.

When writing partitions, we usually group together equal parts with a superscript and omit the trailing zeroes, so that $(4,3,3,1,1,0,0,0,\dots)$ is written as $(4,3^2,1^2)$. The partition $(0,0,0,\dots)$ is written as $\varnothing$.

We will often consider partitions in which a subsequence of the parts is an arithmetic progression with common difference $4$. So given integers $a\gs b$ with $a\equiv b\ppmod4$, we write $a\df b$ for the arithmetic progression $a,a-4,a-8,\dots,b$. For example, $(17,15\df3,2)$ means the partition $(17,15,11,7,3,2)$.

The \emph{Young diagram} of a partition $\la$ is the set
\[
[\la]=\lset{(r,c)\in\bbn^2}{c\ls\la_r}
\]
whose elements are called the \emph{nodes} of $\la$. We say that a node of $\la$ is \emph{removable} if it can be removed from $[\la]$ to leave the Young diagram of a partition; on the other hand, a node not in $[\la]$ is an addable node of $\la$ if it can be added to $\la$ to give the Young diagram of a partition.

We depict Young diagrams as arrays of boxes using the English convention, in which $r$ increases down the page and $c$ increases from left to right. We often blur the distinction between a partition and its Young diagram; for example, we may write $\la\subseteq\mu$ to mean that $[\la]\subseteq[\mu]$.

If $\la$ is a partition, the \emph{conjugate} partition $\la'$ is given by
\[
\la'_r=\card{\lset{c\in\bbn}{\la_c\gs r}}.
\]
In other words, $\la'$ is the partition whose Young diagram is obtained by reflecting $[\la]$ on the main diagonal. For example, we write $(4,3,1^2)'=(4,2^2,1)$.

A partition $\la$ is \emph{\treg}  (or \emph{strict}) if its positive parts are all different, and \emph{$2$-singular} otherwise. We write $\cald$ for the set of \trps, and $\cald(n)$ for the set of \trps of $n$.

Now we recall some additional notation from \cite{mfspin2}, which is very natural but not quite standard. Suppose we have partitions $\la$ and $\mu$ and a natural number $m$. Then we write:
\begin{itemize}
\item
$m\la$ for the partition $(m\la_1,m\la_2,\dots)$;
\item
$\la+\mu$ for the partition $(\la_1+\mu_1,\la_2+\mu_2,\dots)$;
\item
$\la\sqcup\mu$ for the partition obtained by combining all the parts of $\la$ and $\mu$ and arranging them in decreasing order.
\end{itemize}

We may combine these operations, and they take precedence in the order they appear above, so that $\la\sqcup m\mu+\nu$ means $\la\sqcup((m\mu)+\nu)$.

We will need the following simple lemma.

\begin{lemma}\label{cupcont}
Suppose $\la,\mu,\nu\in\calp$. Then $\la\sqcup\mu\subseteq\la\sqcup\nu$ \iff $\mu\subseteq\nu$.
\end{lemma}

\begin{pf}
This is easy to see using conjugate partitions: observe that
\[
(\la\sqcup\mu)'=\la'+\mu',\qquad(\la\sqcup\nu)'=\la'+\nu'.
\]
So
\begin{align*}
\la\sqcup\mu\subseteq\la\sqcup\nu&\Longleftrightarrow(\la\sqcup\mu)'\subseteq(\la\sqcup\nu)'\\
&\Longleftrightarrow\la'_r+\mu'_r\ls\la'_r+\nu'_r\quad\text{for all $r$}\\
&\Longleftrightarrow\mu'_r\ls\nu'_r\quad\text{for all $r$}\\
&\Longleftrightarrow\mu'\subseteq\nu'\\
&\Longleftrightarrow\mu\subseteq\nu.\qedhere
\end{align*}
\end{pf}

The \emph{dominance order} is a partial order defined as follows: given compositions $\la$, $\mu$, we say that $\la$ dominates $\mu$ (and write $\la\dom\mu$) if $|\la|=|\mu|$ and
\[
\la_1+\dots+\la_r\gs\mu_1+\dots+\mu_r
\]
for all $r$.

We shall be mainly interested in the restriction of $\dom$ to $\calp$. The following easy lemma is well known.

\begin{lemma}\label{domcover}
Suppose $\la,\mu\in\calp$. Then $\la$ covers $\mu$ in the dominance order on $\calp$ \iff $\mu$ is obtained from $\la$ by moving one node either to the row immediately below, or to the column immediately to the left.
\end{lemma}

We also need a lemma relating the dominance order and the operator $\sqcup$; the proof is almost identical to the proof of \cref{cupcont}, but using the fact that $\mu\dom\nu$ \iff $\nu'\dom\mu'$.

\begin{lemma}\label{domsqcup}
Suppose $\la,\mu,\nu\in\calp$. Then $\la\sqcup\mu\dom\la\sqcup\nu$ \iff $\mu\dom\nu$.
\end{lemma}

We shall also occasionally need to consider \emph{skew partitions}. A skew partition is a pair of partitions $\la,\mu$ such that $\la\supseteq\mu$. We write this pair as $\la\setminus\mu$. The Young diagram of $\la\setminus\mu$ is then $[\la]\setminus[\mu]$. We identify skew partitions which have the same Young diagram; so for example $(4,2^2)\setminus(3,2,1)=(4,3,2,1)\setminus(3^3,1^2)$. Furthermore, we regard any partition $\la$ as a skew partition by identifying it with $\la\setminus\varnothing$.

\subsection{$2$-cores, residues and $2$-quotients}\label{corequotsec}

Here we recall some of the combinatorics associated with the $2$-modular representation theory of the symmetric group. A \emph{rim $e$-hook} (also called an \emph{$e$-ribbon}) is a set of $e$ nodes $\{(r_1,c_1),\dots,(r_e,c_e)\}$ in $\bbn^2$ such that $(r_{i+1},c_{i+1})$ equals either $(r_i,c_i+1)$ or $(r_i-1,c_i)$, for each $i$. In particular, a rim $2$-hook is just a pair of horizontally or vertically adjacent nodes.

A \emph{rim $2$-hook} of a partition $\la$ is a rim $2$-hook contained in $\la$ which can be removed to leave the Young diagram of a smaller partition. A \emph{$2$-core} is a partition with no rim $2$-hooks, i.e.\ a partition of the form $(c,c-1,\dots,1)$ for some $c\gs0$. If $\la$ is any partition, the $2$-core of $\la$ is the partition obtained by repeatedly removing rim $2$-hooks until none remain.

$2$-cores are closely related to residues. The \emph{residue} of a node $(r,c)$ is the residue of the integer $c-r$ modulo $2$. We call a node of residue $i$ an \emph{$i$-node}. We define the \emph{content} of a partition $\la$ to be the multiset of the residues of the nodes of $\la$. We write the content of a partition in the form $\{0^a,1^b\}$; for example, the partition $(5,2,1)$ has content $\{0^5,1^3\}$. The following lemma was proved by Littlewood \cite{litt}.

\begin{lemma}\label{corecont}
Suppose $\la$ and $\mu$ are partitions of the same size. Then $\la$ and $\mu$ have the same $2$-core \iff they have the same content.
\end{lemma}

We also want to assign a sign to every partition with empty $2$-core. If the $2$-core of $\la$ is $\varnothing$, then we define $\epstar\la=(-1)^h$, where $h$ is the number of horizontal rim $2$-hooks removed in reducing $\la$ to its $2$-core. (In fact $h$ is not well-defined, but its parity is, which is all we need.) Observe that $\epstar\la$ is just the \emph{$2$-sign} (as defined in \cite[p.229]{j10}) of $\la'$. However, in contrast to \cite{j10}, we set $\epstar\la=0$ if the $2$-core of $\la$ is not $\varnothing$.

We will also need the $2$-quotient of a partition, and for this it is convenient to use the abacus, introduced by James \cite{jk}. Take an abacus with two vertical runners labelled $0$ and $1$, and mark positions labelled with the integers on these runners, so that all even integers appear on runner $0$ and all the odd integers appear on runner $1$, increasing from top to bottom. The \emph{abacus display} for a partition $\la$ is obtained by placing a bead on the abacus at position $\la_r-r$ for each $r$. Now for $a\in\{0,1\}$, let $\la^{(a)}_r$ be the number of empty positions above the $r$th lowest bead on runner $a$, for each $r$. Then $\quo\la a=(\quo\la a_1,\quo\la a_2,\dots)$ is a partition, and the pair $(\quo\la0,\quo\la1)$ is the \emph{$2$-quotient} of $\la$. An abacus display for the $2$-core of $\la$ can be obtained by sliding all beads up until every bead has a bead immediately above it.

For example, the partition $(5,3,1^3)$ has $2$-core $(1)$ and $2$-quotient $((2,1),(2))$, as we see from its abacus display.
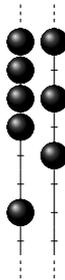
\[
\abacus(vv,bb,bn,bb,bn,nb,nn,bn,nn,vv)
\]

We will later need the following two lemmas.  These are probably not new, but the author has not been able to find them in the literature.

\begin{lemma}\label{duplem}
Suppose $\alpha$ is a partition and $c\gs1$ with $\alpha'_c>\alpha'_{c+1}$. Let $\beta=\dup\alpha$ and let $\gamma$ be the partition obtained from $\beta$ by moving a node from column $c$ to column $c+1$. Then $\alpha=\quo\beta0+\quo\beta1=\quo\gamma0+\quo\gamma1$, and $\epstar\gamma=-1$.
\end{lemma}

\begin{pf}
First we show that $\quo\beta0+\quo\beta1=\alpha$. Since $\beta_{2r-1}=\beta_{2r}$ for each $r$, the $(2r-1)$th and $2r$th lowest beads in the abacus display for $\beta$ occur in consecutive positions, namely positions $\alpha_r-2r,\alpha_r-2r+1$. So these two beads are (in some order) the $r$th lowest bead on runner $0$ and the $r$th lowest bead on runner $1$. So (by the definition of $2$-quotient) $\quo\beta0_r+\quo\beta1_r$ equals the total number of empty positions before position $\alpha_r-2r$, which is $\alpha_r$. So $\quo\beta0+\quo\beta1=\alpha$, as claimed.

Now consider the abacus display for $\gamma$. This is obtained from the abacus display for $\beta$ by moving two beads: writing $s=\alpha'_{c+1}$ and $t=\alpha'_c$, the abacus display for $\beta$ includes beads at positions
\begin{align*}
\beta_{2s+1}=c-2s-1,\,\beta_{2s+2}=c-2s-2&,\,\dots,\,\beta_{2t-1}=c-2t+1,\,\beta_{2t}=c-2t
\\
\intertext{and to construct $\gamma$ these are replaced by beads at positions}
\gamma_{2s+1}=c-2s,\ \gamma_{2s+2}=c-2s-2&,\dots,\,\gamma_{2t-1}=c-2t+1,\ \gamma_{2t}=c-2t-1.
\end{align*}
Looking at the positions of each parity in turn, this yields
\[
\quo\gamma i_r=
\begin{cases}
\quo\beta i_r+1&(r\in\{s+1,\dots,t\},\ c\equiv i\ppmod2)\\
\quo\beta i_r-1&(r\in\{s+1,\dots,t\},\ c\not\equiv i\ppmod2)\\
\quo\beta i_r&(\text{otherwise}).
\end{cases}
\]
Hence $\quo\gamma0+\quo\gamma1=\quo\beta0+\quo\beta1$.

We show that $\epstar\gamma=-1$ by showing that $\gamma$ can be reduced to the empty partition by removing rim $2$-hooks, with exactly one of the removed rim $2$-hooks being horizontal.  The construction of $\gamma$ means that $\gamma'_c-\gamma'_{c+1}=2(t-s-1)$. So we begin by removing $t-s-1$ vertical rim $2$-hooks from column $c$. We can then remove a horizontal rim $2$-hook $\{(2s+1,c),(2s+1,c+1)\}$, to leave a partition in which all the columns have even length. Such a partition can clearly be reduced to $\varnothing$ by repeatedly removing vertical rim $2$-hooks.
\end{pf}

\begin{lemma}\label{empcore}
Suppose $\gamma$ is a partition with $2$-core $\varnothing$. Then $\gamma\dom(\quo\gamma0+\quo\gamma1)\sqcup(\quo\gamma0+\quo\gamma1)$, with equality \iff all the columns of $\gamma$ have even length.
\end{lemma}

\begin{pf}
In the abacus display for $\gamma$, let $b_r^i$ be the position of the $r$th lowest bead on runner $i$. If $|b^0_r-b^1_r|=1$ for all $r$, then $\gamma_{2r-1}=\gamma_{2r}$ for all $r$, so that all the columns of $\gamma$ have even length. This means that we can write $\gamma$ in the form $\dup\alpha$, and \cref{duplem} (with $\gamma$ in place of $\beta$) shows that $\alpha=\quo\gamma0+\quo\gamma1$, as claimed.

Now suppose the columns of $\gamma$ do not all have even length. Then there is some $r$ for which $|b^0_r-b^1_r|>1$; let $k$ be the maximum value of $|b^0_r-b^1_r|$, and take $s$ minimal such that $|b^0_s-b^1_s|=k$. We assume for simplicity that $b^0_s-b^1_s=k$; the other case is identical, but with $0$ and $1$ interchanged.

The choice of $s$ means that position $b_s^1+2$ in the abacus is empty. Now let $d\gs0$ be maximal such that positions $b^0_s-2,\dots,b^0_s-2d$ are all occupied. The maximality of $k$ then means that positions $b^1_s-2,\dots,b^1_s-2d$ are also all occupied. We construct a new partition $\eta$ by moving the bead at position $b^0_s$ to position $b^0_s-2d-2$, and moving the bead at position $b^1_s-2d$ to position $b^1_s+2$. This corresponds to removing a rim $(2d+2)$-hook from $\gamma$, and then adding a rim $(2d+2)$-hook in a lower position (see \cite[\S2.7]{jk} for details of rim hooks and the abacus), so that $\eta\domsby\gamma$. We also have $\quo\eta0+\quo\eta1=\quo\gamma0+\quo\gamma1$; this is very similar to the calculation of $\quo\gamma0+\quo\gamma1$ in the proof of \cref{duplem}. By induction on the dominance order $\eta\dom\dup{(\quo\eta0+\quo\eta1)}$, and so $\gamma\doms\dup{(\quo\gamma0+\quo\gamma1)}$.
\end{pf}

\subsection{Regularisation, doubling and \fbcs}\label{regdoubsec}

Here we introduce two operations on partitions which have significance for decomposition numbers of symmetric groups in characteristic $2$.

For $l\gs0$, we define the $l$th \emph{ladder} in $\bbn^2$ to be the set of nodes $(r,c)$ for which $r+c=l+2$. Given a partition $\la$, its \emph{regularisation} $\la\reg$ is the \trp obtained by moving the nodes of $\la$ as high as they will go within each ladder. For example, the regularisation of $(3^2,1^3)$ is $(5,3,1)$, as we see from the following diagrams, in which we label each node with the number of the ladder in which it lies.
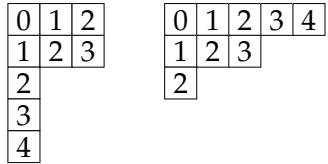
\[
\young(012,123,2,3,4)\qquad\young(01234,123,2)
\]

If $\la$ is a \trp, we define its \emph{double} $\la\dbl$ to be the partition
\[
(\lfloor\la_1/2\rfloor,\lceil\la_1/2\rceil,\lfloor\la_2/2\rfloor,\lceil\la_2/2\rceil,\dots).
\]
Regularisation was introduced by James \cite{j1}, and doubling by Bessenrodt and Olsson \cite{bo}. We will see the representation-theoretic significance of these operations in \cref{rdsec}. For a \trp $\la$ we write $\la\dblreg$ for $(\la\dbl)\reg$.

A \emph{\fbc} is a partition of the form $(4l-1\df3)$ or $(4l-3\df1)$ for some $l\gs0$. Observe that a partition is a \fbc precisely if its double is a $2$-core (or equivalently, if its conjugate is the double of a $2$-core). Given a \trp $\la$, we define the \fbc of $\la$ to be the \fbc obtained by repeatedly applying the following operations:
\begin{itemize}
\item
removing even parts from $\la$;
\item
removing any two parts whose sum is a multiple of $4$;
\item
replacing any odd part $\la_i\gs5$ with $\la_i-4$, if $\la_i-4$ is not already a part of $\la$.
\end{itemize}

Now we have the following.

\begin{lemma}\label{regdoub}
Suppose $\la$ and $\mu$ are \trps and $\tau$ is the \fbc of $\la$. Then the $2$-core of $\la\dblreg$ is $\tau\dbl$. Hence $\la$ and $\mu$ have the same \fbc \iff $\la\dblreg$ and $\mu\dblreg$ have the same $2$-core.
\end{lemma}

\begin{pf}
By \cite[Lemma 3.6]{bo} the $2$-core of $\la\dbl$ is $\tau\dbl$. But $\la\dbl$ and $\la\dblreg$ have the same content, because all the nodes in a given ladder have the same residue. Hence by \cref{corecont} $\la\dbl$ and $\la\dblreg$ have the same $2$-core, which gives the first statement. The second statement then follows from the statement that the function $\tau\mapsto\tau\dbl$ yields a bijection from the set of \fbcs to the set of $2$-cores: this bijection is easily seen just by writing it down directly:
\begin{alignat*}2
\varnothing&\longmapsto\varnothing&(1)&\longmapsto(1)
\\
(3)&\longmapsto(2,1)&(5,1)&\longmapsto(3,2,1)
\\
(7,3)&\longmapsto(4,3,2,1)&\qquad(9,5,1)&\longmapsto(5,4,3,2,1)
\\
&\ \ \ \vdots&&\ \ \ \vdots\tag*{\qedhere}
\end{alignat*}
\end{pf}

\subsection{Semistandard tableaux}\label{semistsec}

Here we give some basic background on tableaux. These are treated at length elsewhere, so we give the minimum amount of detail required for this paper. The book by Fulton \cite{fult} is an excellent reference.

Suppose $\la$ is a skew partition and $\Omega$ is a set. A \emph{$\la$-tableau} with entries in $\Omega$ is a function $T:[\la]\to\Omega$. We think of a tableau $T$ as a way of filling the boxes of $[\la]$ with elements of $\Omega$, and we say that $T$ has \emph{shape} $\la$.

Now suppose $\Omega$ is equipped with a total order $\ls$. We say that a $\la$-tableau $T$ is \emph{semistandard} if the entries in each row are weakly increasing (with respect to $\ls$) from left to right and the entries in each column are strictly increasing from top to bottom.

If the totally ordered set $\Omega$ is not specified, it should be taken to be $\bbn$ with the usual ordering. In this case, given a tableau $T$, we define a composition $\mu$ by letting $\mu_r$ equal the number of entries equal to $r$ in $T$ for each $r$; then we say that $T$ has \emph{type} $\mu$. In the case where $\la$ is a partition, we write $K_{\la\mu}$ for the number of semistandard $\la$-tableaux of type $\mu$. Then $K_{\la\mu}$ is called a \emph{Kostka number}. It is easy to see that $K_{\la\la}=1$, while $K_{\la\mu}=0$ unless $\la\dom\mu$.

\subsection{Littlewood--Richardson coefficients}\label{lrsec}

Here we define Littlewood--Richardson coefficients. Suppose $\la$ is a skew partition and $T$ is a $\la$-tableau (with entries in $\bbn$). The \emph{reading word} of $T$ is the word obtained by reading the entries \emph{from right to left} along successive rows of the tableau from top to bottom. We say that $T$ is \emph{Littlewood--Richardson} if it is semistandard and its reading word $w$ is a \emph{lattice word}, meaning that every initial segment of $w$ contains at least as many $r$s as $(r+1)$s, for every~$r$.

For example, with $\la=(4,3,2)\setminus(2)$, the following are all semistandard tableaux of type $(3,2^2)$, but only the first is Littlewood--Richardson. We give the reading word in each case.
\[
\begin{array}{c@{\qquad}c@{\qquad}c}
\gyoung(::;11,122,33)&\gyoung(::;11,123,23)&\gyoung(::;12,112,33)\\[36pt]
1122133&1132132&2121133
\end{array}
\]

Now suppose $\alpha,\beta,\gamma$ are partitions. The \emph{Littlewood--Richardson coefficient} $\lrc\gamma\alpha\beta$ is defined to be the number of Littlewood--Richardson tableaux of shape $\gamma\setminus\alpha$ and type $\beta$ if $\alpha\subseteq\gamma$, and $0$ otherwise. Obviously $\lrc\gamma\alpha\beta=0$ unless $|\gamma|=|\alpha|+|\beta|$.

Two special cases of the Littlewood--Richardson coefficients arise where the partition $\beta$ has only one row or one column. Say that a skew partition $\gamma\setminus\alpha$ is a \emph{horizontal $a$-strip} if $[\gamma]$ is obtained from $[\alpha]$ by adding $a$ nodes in distinct columns, or a \emph{vertical $a$-strip} if $[\gamma]$ is obtained from $[\alpha]$ by adding $a$ nodes in distinct rows. Then
\begin{align*}
\lrc\gamma\alpha{(a)}&=
\begin{cases}
1&\text{if $\gamma\setminus\alpha$ is a horizontal $a$-strip}\\
0&\text{otherwise},
\end{cases}\\
\lrc\gamma\alpha{(1^a)}&=
\begin{cases}
1&\text{if $\gamma\setminus\alpha$ is a vertical $a$-strip}\\
0&\text{otherwise}.
\end{cases}
\end{align*}
These  formul\ae{} are called the \emph{Pieri rule} and the \emph{dual Pieri rule}. Later we shall use the following notation: given partitions $\alpha,\gamma$, we write $\lrcb\gamma\alpha$ for $\sum_{a\gs0}\lrc\gamma\alpha{(1^a)}$. In other words, $\lrcb\gamma\alpha=1$ if $\gamma\supseteq\alpha$ and $\gamma\setminus\alpha$ is a vertical strip of any size, and otherwise $\lrcb\gamma\alpha=0$.

We will need the following basic results concerning Littlewood--Richardson coefficients.

\begin{lemmac}{j10}{(2.3)}\label{basiclr}
Suppose $\alpha,\beta,\gamma\in\calp$. Then
\[
\lrc{\alpha+\beta}\alpha\beta=\lrc{\alpha\sqcup\beta}\alpha\beta=1
\]
while $\lrc\gamma\alpha\beta=0$ unless $\alpha+\beta\dom\gamma\dom\alpha\sqcup\beta$.
\end{lemmac}

\begin{lemma}\label{lrlem}
Suppose $\la,\beta,\zeta\in\calp$ and $r\gs1$ such that $\la_r>\la_{r+1}$, and $\zeta$ is obtained from $\beta$ by moving a node from row $r$ to row $r+1$. Then $\lrc{\la+\zeta}\la\beta>0$.
\end{lemma}

\begin{pf}
For this, we need to construct a Littlewood--Richardson tableau of shape $(\la+\zeta)\setminus\la$ and content $\beta$. Let $T$ be the unique Littlewood--Richardson tableau $T$ of shape $(\la+\zeta)\setminus\la$ and content $\zeta$; this has all entries in row $s$ equal to $s$, for each $s$. Now let $U$ be the tableau obtained by replacing the first entry in row $r+1$ with $r$. Then $U$ has content $\beta$, and is semistandard because $\la_r>\la_{r+1}$. Furthermore, it is easy to check that $U$ is a Littlewood--Richardson tableau, so $\lrc{\la+\zeta}\la\beta>0$.
\end{pf}

\section{Symmetric functions}\label{symfnsec}

In this section we introduce the background we need on symmetric functions, and prove an auxiliary result that we shall need for the proof of our theorem on decomposition numbers for Rouquier blocks. The standard reference for symmetric functions is Macdonald's book \cite{macd}; we recall here the basics needed for this paper.

\subsection{The algebra of symmetric functions}\label{algsymfnsec}

We let $X$ be a countably infinite set of commuting algebraically independent indeterminates, and we define $\La$ to be the $\bbq$-algebra of symmetric functions in $X$; that is, power series of bounded degree which are invariant under all permutations of $X$. For any countably infinite set $Y$ of indeterminates and any $f\in\La$, we can define $f(Y)$ by replacing the elements of $X$ with the elements of $Y$.

$\La$ is equipped with a coproduct $\Delta:\La\to\La\otimes\La$, defined as follows. We partition $X$ into two infinite sets $X=Y\sqcup Z$. Then $f\in\La$ is a symmetric function of the elements of $Y$ and also a symmetric function of the elements of $Z$. So we can write $f$ as a finite sum $f=\sum_i f_i(Y)g_i(Z)$ where $f_i,g_i\in\La$ for all $i$. We now define $\Delta(f)=\sum_i f_i\otimes g_i$.

\subsubsection*{Monomial functions and Schur functions}

$\La$ is equipped with several important bases, and we recall some of the details here. We begin with the \emph{monomial symmetric functions}. Given a partition $\la$ of length $l$, define $m_\la$ to be the sum of all distinct monomials of the form $x_1^{\la_1}\dots x_l^{\la_l}$, where $x_1,\dots,x_l$ are distinct elements of $X$. Clearly $m_\la$ is a symmetric function which is homogeneous of degree $|\la|$, and the set $\lset{m_\la}{\la\in\calp}$ is a basis for $\La$.

We also recall the definition of the \emph{Schur functions}. Suppose $\la\in\calp$, and choose a total order $\ls$ on $X$. We write $\sstd\la$ for the set of all semistandard $\la$-tableaux with entries in $X$. If $T\in\sstd\la$, we define the monomial $x^T$ to be the product of the entries in $T$. The Schur function $s_\la$ is defined to be
\[
s_\la=\sum_{T\in\sstd\la}x^T.
\]
Then $s_\la$ is a symmetric function which is homogeneous of degree $|\la|$ and independent of the choice of total order on $X$.

The Schur functions comprise a basis of $\La$. So we can define an inner product $\ip\ \ $ on $\La$ by specifying that the Schur functions are orthonormal. There is also a unique linear map $\omega:\La\to\La$ satisfying $\omega(s_\la)=s_{\la'}$ for all $\la$. In fact, $\omega$ is an algebra automorphism \cite[I.2]{macd}.

The transition between Schur and monomial symmetric functions is well understood: we have $s_\la=\sum_{\mu\in\calp}K_{\la\mu}m_\mu$, where $K_{\la\mu}$ is the Kostka number introduced in \cref{semistsec}.

By \cite[I.9.2]{macd}, the structure constants for the basis of Schur functions are the Littlewood--Richardson coefficients; that is,
\[
s_\alpha s_\beta=\sum_{\gamma\in\calp}\lrc\gamma\alpha\beta s_\gamma
\]
for all $\alpha,\beta\in\calp$. As a consequence, we have $\lrc\gamma\alpha\beta=\lrc\gamma\beta\alpha$ for all $\alpha,\beta,\gamma$. The fact that $\omega$ is an automorphism also gives $\lrc\gamma\alpha\beta=\lrc{\gamma'}{\alpha'}{\beta'}$.

Given this interpretation of Littlewood--Richardson coefficients in terms of Schur functions, it is natural to extend the definition of these coefficients to allow more than three arguments: we may define $\lrc\zeta{\alpha\beta}\gamma$ to be the coefficient of $s_\zeta$ in the product $s_\alpha s_\beta s_\gamma$, or equivalently $\lrc\zeta\alpha{\beta\gamma}=\sum_{\mu\in\calp}\lrc\zeta\alpha\mu\lrc\mu\beta\gamma$. We also define $\lrcb\zeta{\alpha\beta}=\sum_{a\gs0}\lrc\zeta{\alpha\beta}{(1^a)}$.

The effect of the coproduct $\Delta$ on Schur functions can also be described in terms of Littlewood--Richardson coefficients: we have
\[
\Delta(s_\la)=\sum_{\mu,\nu\in\calp}\lrc\la\mu\nu s_\mu\otimes s_\nu
\]
for any $\la\in\calp$ \cite[I.5, Example 25]{macd}. A consequence of this (and the conjugation-symmetry of the Littlewood--Richardson coefficients) is that if $f\in\La$ with $\Delta(f)=\sum_if_i\otimes g_i$, then $\Delta(\omega(f))=\sum_i\omega(f_i)\otimes\omega(g_i)$.

\subsubsection*{Schur P-functions}

Next we recall the family of \emph{Schur P-functions}. Our definition (modulo a minor change of notation) follows Stembridge \cite[\S6]{stem}.

Suppose $\la$ is a \trp. The \emph{shifted Young diagram} of $\la$ is the set
\[
\shyng\la=\lset{(r,c)\in\bbn^2}{r\ls c<r+\la_r}.
\]
As above, we choose a total order  $\ls$ on $X$. We define $X^\pm$ to be the set of symbols $\lset{x^+,x^-}{x\in X}$, and we define a total order (also called $\ls$) on $X^\pm$ by setting $x^-<x^+$ for all $x\in X$, and $x^\pm<y^\pm$ for all $x,y\in X$ with $x<y$.

Now we define a \emph{semistandard shifted $\la$-tableau} to be a filling of the boxes of $\shyng\la$ with elements of $X^\pm$ such that:
\begin{itemize}
\item
the symbols are weakly increasing (with respect to $\ls$) down each column and from left to right along each row;
\item
for each $x\in X$, the symbol $x^-$ appears at most once in each row, and the symbol $x^+$ appears at most once in each column;
\item
the entries on the diagonal are all of the form $x^+$.
\end{itemize}
We write $\shtd\la$ for the set of semistandard shifted $\la$-tableaux. Given $T\in\shtd\la$, we define the monomial $x^T$ to be the product of the symbols in $T$ (with the signs $\pm$ ignored). The Schur P-function labelled by $\la$ is then defined as
\[
P_\la=\sum_{T\in\shtd\la}x^T.
\]
By \cite[Corollary 6.2(a)]{stem} $P_\la$ is a symmetric function which is homogeneous of degree $|\la|$.

We shall need the following results, which are surprisingly difficult to find explicitly in the literature.

\begin{lemma}\label{pomeg}
Suppose $\la\in\cald$. Then $\omega(P_\la)=P_\la$.
\end{lemma}

\begin{pf}
For $r\gs0$, let $p_r$ denote the \emph{power-sum symmetric function} $\sum_{x\in X}x^r$. From \cite[p.24]{macd} we have $\omega(p_r)=p_r$ whenever $r$ is odd. Hence the algebra $\Omega$ generated by $p_1,p_3,p_5,\dots$ consists entirely of $\omega$-invariant symmetric functions. By \cite[Corollary 6.2(b)]{stem} the Schur P-functions comprise a basis for $\Omega$.
\end{pf}

\begin{lemma}\label{schurpdom}
Suppose $\la\in\cald$ and $\mu\in\calp$. Then $\lan P_\la,s_\la\ran=1$, and $\lan P_\la,s_\mu\ran=0$ unless $\la\dom\mu$.
\end{lemma}

\begin{pf}
By \cite[Corollary 6.2 and Lemma 6.3]{stem} the transition matrix from Schur P-functions to monomial symmetric functions is unitriangular, in the sense that $P_\la$ equals $m_\la$ plus a linear combination of functions $m_\mu$ with $\mu\domsby\la$. The definition of Schur functions implies that the transition matrix from monomial symmetric functions to Schur functions is also unitriangular, and hence the transition matrix from Schur P-functions to Schur functions is unitriangular.
\end{pf}

\begin{cory}\label{p2core}
Suppose $\la$ is a $2$-core. Then $P_\la=s_\la$.
\end{cory}

\begin{pf}
By \cref{pomeg} and the definition of $\omega$, we get $\lan P_\la,s_\mu\ran=\lan P_\la,s_{\mu'}\ran$ for every $\mu$. So by \cref{schurpdom} we must have $\la\dom\mu$ and $\la\dom\mu'$ in order to get $\lan P_\la,s_\mu\ran\neq0$, or in other words $\la\dom\mu\dom\la'$. In the case where $\la$ is a $2$-core, this gives $\la\dom\mu\dom\la$, so that $\mu=\la$.
\end{pf}

\begin{lemma}\label{schurpcover}
Suppose $\la\in\cald$ and $r\gs1$ such that $\la_r\gs\la_{r+1}+2$, and let $\mu$ be the partition obtained from $\la$ by moving a node from row $r$ to row $r+1$. Then $\lan P_\la,s_\mu\ran=1$.
\end{lemma}

\begin{pf}
First we consider the expansion $P_\la=\sum_\nu C_{\la\nu}m_\nu$ as a sum of monomial symmetric functions and find the coefficient $C_{\la\mu}$. To do this, we write $X=\{x_1,x_2,\dots\}$ with $x_1<x_2<\dots$, and we just need to find all semistandard shifted $\la$-tableaux $T$ such that $x^T=x_1^{\mu_1}x_2^{\mu_2}\dots$. It is easy to see that there are exactly two such tableaux, namely the tableaux in which all the entries in row $s$ are equal to $x_s^+$ for every $s$, with the exception of the last entry in row $r$, which is either $x_{r+1}^+$ or $x_{r+1}^-$. (Both possibilities can occur, since $\la_r\gs\la_{r+1}+2$.)

Hence $C_{\la\mu}=2$. Since $\la$ covers $\mu$ in the dominance order and $C_{\la\nu}=0$ when $\la\ndom\nu$, we can write
\[
P_\la=m_\la+2m_\mu+\sum_{\nu\ndom\mu}C_{\la\nu}m_\nu.
\]

In a similar way there is exactly one semistandard $\la$-tableau $T$ with $x_T=x_1^{\mu_1}x_2^{\mu_2}\dots$, so that
\[
s_\la=m_\la+m_\mu+\sum_{\nu\ndom\mu}K_{\la\nu}m_\nu.
\]
Hence $s_\mu$ occurs exactly once in $P_\la$.
\end{pf}

\subsection{A reduction operator on symmetric functions}\label{redopsec}

For any partition $\nu$ we define a linear map $\ddd_\nu$ on $\La$ by mapping
\[
s_\la\mapsto\sum_{\mu\in\calp}\lrc\la\mu\nu s_\mu
\]
and extending linearly. Then $\ddd_\nu$ reduces the degree of a homogeneous function by $|\nu|$. We shall be particularly interested in the cases $\nu=(a)$ and $\nu=(1^a)$. We will use the following lemma.

\begin{lemma}\label{smd}
Suppose $f\in\La$, and write
\[
\Delta(f)=\sum_{\mu\in\calp}g_\mu\otimes m_\mu,
\]
where $g_\mu\in\La$ for each $\mu$. If $a$ is any non-negative integer, then $\ddd_{(a)}(f)=g_{(a)}$.
\end{lemma}

\begin{pf}
Recalling that $\Delta(s_\la)=\sum_{\mu,\nu\in\calp}\lrc\la\mu\nu s_\mu\otimes s_\nu$ for every $\la$ and using the definition of $\ddd_\nu$, we obtain
\[
\Delta(f)=\sum_{\nu\in\calp}\ddd_\nu(f)\otimes s_\nu.
\]
Writing $s_\nu=\sum_\mu K_{\nu\mu}m_\mu$, we have
\[
\Delta(f)=\sum_{\mu,\nu}K_{\nu\mu}\ddd_\nu(f)\otimes m_\mu.
\]
Since the $m_\mu$ are linearly independent, this gives
\[
g_\mu=\sum_\nu K_{\nu\mu}\ddd_\nu(f)
\]
for each $\mu$, and in particular $g_{(a)}=\ddd_{(a)}(f)$.
\end{pf}

As a consequence, it is easy to determine the effect of the function $\ddd_{(a)}$ on a symmetric function $f$. \cref{smd} says that  is the function $g_{(a)}$ appearing in the expansion $\Delta(f)=\sum_{\mu\in\calp}g_\mu\otimes m_\mu$. To find this expansion, partition $X$ as the union $Y\sqcup Z$ of two infinite sets; now write each monomial in $f$ as the product of a monomial function of the elements of $Y$ and a monomial function of the elements of $Z$. For given $\mu$, taking the sum of all the terms for which the ``$Z$-part'' has the form $z_1^{\mu_1}z_2^{\mu_2}\dots$ for distinct $z_1,z_2,\dots$, we obtain the summand $g_\mu\otimes m_\mu$. In particular, for the case $\mu=(a)$, we take all monomials in $f$ in which the $Z$-part equals $z^a$ for some $z\in Z$, and the sum of these monomials is $g_{(a)}\otimes m_{(a)}=g_{(a)}\otimes\sum_{z\in Z}z^a$. So to write down $g_{(a)}$, we can just fix an element $z\in Z$, and find all monomials in $f$ of the form $z^a$ times a monomial in $Y$. The sum of these monomials is then $z^ag_{(a)}(Y)=z^a(\ddd_{(a)}f)(Y)$.

We want to describe the effect of the linear maps $\ddd_{(a)}$ and $\ddd_{(1^a)}$ on Schur P-functions. To this end, we introduce some notation. Suppose $\la$ and $\mu$ are \trps. If $\la\setminus\mu$ is a horizontal strip, then we define $N_{\la\mu}$ to be the number of $c\gs1$ such that $\la\setminus\mu$ contains a node in column $c+1$ but not in column $c$, and set $\strp\la\mu=2^{N_{\la\mu}}$. Otherwise, we set $\strp\la\mu=0$. (Our choice of notation reflects the fact that $\strp\la\mu$ is similar to the quantity $\operatorname{hs}_{\la\setminus\mu}(-1)$ introduced by Konvalinka and Lauve \cite{konlau}, but slightly different: if $\len\la>\len\mu$, then $\operatorname{hs}_{\la\setminus\mu}(-1)=2\strp\la\mu$.)

Now we have the following result. This is a kind of Pieri rule for Schur P-functions which appears (surprisingly) to be new. (A related rule appears in \cite[III.5.7]{macd}, but we cannot see how to deduce our result from this.)

\begin{propn}\label{drp}
Suppose $\la\in\cald(n)$ and $a\gs0$. Then
\[
\ddd_{(a)}P_\la=\ddd_{(1^a)}P_\la=\sum_{\mu\in\cald(n-a)}\strp\la\mu P_\mu.
\]
\end{propn}

\begin{pf}
First we consider $\ddd_{(a)}P_\la$. We calculate this as explained above: we write $X$ as the union of two disjoint infinite sets $Y\sqcup Z$, fix $z\in Z$, and consider all the monomials in $P_\la$ of the form $z^a$ times a monomial in $Y$. When we define $P_\la$ in terms of shifted semistandard tableaux, we can choose the total order $\ls$ on $X$ freely; for this proof we choose $\ls$ in such a way that $y<z$ for all $y\in Y$. Now we need to classify shifted semistandard tableaux $T$ which have $a$ entries of the form $z^\pm$ with the remaining entries lying in $Y^\pm$. Since $y<z$ for all $y\in Y$, the entries in $Y^\pm$ in such a tableau comprise a shifted semistandard $\mu$-tableau for some $\mu\in\cald$ with $\mu\subseteq\la$. The definition of a shifted semistandard tableau means that we cannot have boxes $(r,c)$ and $(r+1,c+1)$ both containing entries $z^\pm$, so the shape defined by the entries $z^\pm$ (that is, the difference $\shyng\la\setminus\shyng\mu$) is a union of non-adjacent rim hooks. This is the same as saying that $\la\setminus\mu$ is a horizontal $a$-strip; moreover, the number of rim hooks comprising $\shyng\la\setminus\shyng\mu$ which do not meet the main diagonal equals the integer $N_{\la\mu}$ defined above.

Conversely, given $\mu\in\cald$ such that $\la\setminus\mu$ is a horizontal $a$-strip, and given a shifted semistandard $\mu$-tableau with entries in $Y^\pm$, we can add $a$ symbols $z^\pm$ to make a shifted semistandard $\la$-tableau. The positions of these symbols are determined by $\la$, and the signs on these symbols are determined by the rules for a shifted semistandard tableau, except for the bottom-left entry in each constituent rim hook of $\shyng\la\setminus\shyng\mu$ not meeting the diagonal, whose sign can be chosen freely.

As a consequence, we see that the sum of the monomials $x^T$, where $T$ contains $a$ symbols $z^\pm$ with the remaining symbols coming from $Y^\pm$, is $\sum_{\mu\in\cald}z^a\strp\la\mu P_\mu(Y)$. So we get $\ddd_{(a)}P_\la=\sum_{\mu\in\cald(n-a)}\strp\la\mu P_\mu$.

To calculate $\ddd_{(1^a)}P_\la$, we use the automorphism $\omega$: the definition of the functions $\ddd_\nu$ and the formula for $\Delta(s_\la)$ means that
\[
\Delta(f)=\sum_{\nu\in\calp}\ddd_\nu f\otimes s_\nu
\]
for any $f\in\La$. Hence
\[
\Delta(\omega(f))=\sum_{\nu\in\calp}\omega(\ddd_\nu f)\otimes s_{\nu'}
\]
and so
\[
\ddd_{\nu'}f=\omega(\ddd_\nu\omega(f)).
\]
We apply this with $\nu=(a)$ and $f=P_\la$. Using \cref{pomeg}, we get
\[
\ddd_{(1^a)}P_\la=\omega(\ddd_{(a)}P_\la)
=\omega\left(\sum_{\mu\in\cald(n-a)}\strp\la\mu P_\mu\right)
=\sum_{\mu\in\cald(n-a)}\strp\la\mu P_\mu.\qedhere
\]
\end{pf}

Now we can derive a corollary which will be an initial case of our main theorem on symmetric functions. Recall that we write $\lrcb\mu\la$ to mean $\sum_{a\gs0}\lrc\mu\la{(1^a)}$.

\begin{cory}\label{specialmain}
Suppose $\alpha\in\cald$ and $\la\in\calp$. Then
\[
\sum_{\mu\in\calp}\lrcb\mu\la\ip{P_\alpha}{s_\mu}=\sum_{\beta\in\cald}\strp\alpha\beta\ip{P_\beta}{s_\la}.
\]
\end{cory}

\begin{pf}
Let $d=|\alpha|$ and $l=|\la|$. If $d<l$ then clearly both sides are zero. Otherwise, let $a=d-l$. Then
\begin{align*}
\sum_{\mu\in\calp}\lrcb\mu\la\ip{P_\alpha}{s_\mu}&=\sum_{\mu\in\calp(d)}\ip{P_\alpha}{s_\mu}\lan\ddd_{(1^a)}s_\mu,s_\la\ran\\
&=\left\lan\ddd_{(1^a)}\sum_{\mu\in\calp(d)}\ip{P_\alpha}{s_\mu}s_\mu,s_\la\right\ran\\
&=\left\lan\ddd_{(1^a)}P_\alpha,s_\la\right\ran\\
&=\left\lan\sum_{\beta\in\cald(l)}\strp\alpha\beta P_\beta,s_\la\right\ran\tag*{by \cref{drp}}\\
&=\sum_{\beta\in\cald}\strp\alpha\beta \lan P_\beta,s_\la\ran.\tag*{\qedhere}
\end{align*}
\end{pf}

\subsection{Main result on symmetric functions}

In this section we will prove our main result on symmetric functions, which we will later use to derive results on decomposition numbers. First we give another initial case for our main theorem. We will deduce this result later from the decomposition number results in \cite{mfspin2}, though it would be preferable to have a proof purely in the context of symmetric functions. Recall the definition of $2$-quotient $(\quo\mu0,\quo\mu1)$ and sign $\epstar\mu$ from \cref{corequotsec}.

\begin{propn}\label{special2main}
Suppose $\la,\gamma\in\calp$. Then
\[
\sum_{\mu\in\calp}\lrcb\mu\la\epstar\mu\lrc\gamma{\quo\mu0}{\quo\mu1}=
\epstar\la\lrcb\gamma{\quo\la0\quo\la1}.
\]
\end{propn}

Now we can give the main theorem of this section.

\begin{thm}\label{mainsym}
Suppose $\la,\gamma\in\calp$ and $\alpha\in\cald$. Then
\[
\sum_{\mu,\nu,\zeta\in\calp}\lrcb\mu\la\ip{P_\alpha}{s_\zeta}\epstar\nu\lrc\gamma{\quo\nu0}{\quo\nu1}\lrc\mu\nu\zeta=
\sum_{\substack{\beta\in\cald\\\eta,\xi\in\calp}}\strp\alpha\beta\ip{P_\beta}{s_\eta}\epstar\xi\lrcb\gamma{\quo\xi0\quo\xi1}\lrc\la\xi\eta.
\]
\end{thm}

We observe that \cref{specialmain} is the case $\gamma=\varnothing$ of \cref{mainsym}. To see this, consider the formula in \cref{mainsym} with $\gamma=\varnothing$. In order to get a non-zero summand on the left-hand side, we must have $\nu=\varnothing$. This then means that $\zeta=\mu$ (otherwise $\lrc\mu\varnothing\zeta=0$) and the left-hand side of the equation in \cref{mainsym} reduces to the left-hand side of the formula in \cref{specialmain}. In order to get a non-zero summand on the right-hand side of the formula in \cref{mainsym} when $\gamma=\varnothing$, we need $\xi=\varnothing$ and hence $\eta=\la$, and we are left with the right-hand side of the formula in \cref{specialmain}. Similarly, we see that \cref{special2main} is the special case $\alpha=\varnothing$ of \cref{mainsym}.

We will prove \cref{mainsym} using \cref{specialmain,special2main}, but first we need one more result about Littlewood--Richardson coefficients; this comes from Mackey's formula applied to characters of symmetric groups.

\begin{lemmac}{ct}{Lemma 2.2(3)}\label{ctmackey}
Suppose $\la,\nu,\zeta\in\calp$. Then
\[
\sum_{\mu\in\calp}\lrcb\mu\la\lrc\mu\nu\zeta=\sum_{\eta,\xi\in\calp}\lrc\la\eta\xi\lrcb\nu\xi\lrcb\zeta\eta.
\]
\end{lemmac}

Now we can prove our main result.

\begin{pf}[Proof of \cref{mainsym}]
Using \cref{ctmackey}, the left-hand side of the formula in \cref{mainsym} equals
\[
\sum_{\eta,\xi,\nu,\zeta\in\calp}\ip{P_\alpha}{s_\zeta}\epstar\nu\lrc\gamma{\quo\nu0}{\quo\nu1}\lrc\la\eta\xi\lrcb\nu\xi\lrcb\zeta\eta.
\]
Using \cref{special2main} with $\nu,\xi$ in place of $\mu,\la$, this becomes
\[
\sum_{\eta,\xi,\zeta\in\calp}\ip{P_\alpha}{s_\zeta}\lrc\la\eta\xi\lrcb\zeta\eta\epstar\xi\lrcb\gamma{\quo\xi0\quo\xi1}.
\]
Using \cref{specialmain} with $\eta,\zeta$ in place of $\la,\mu$, this becomes
\[
\sum_{\substack{\eta,\xi\in\calp\\\beta\in\cald}}\lrc\la\eta\xi\epstar\xi\lrcb\gamma{\quo\xi0\quo\xi1}\strp\alpha\beta\ip{P_\beta}{s_\eta},
\]
which is what we want.
\end{pf}

\section{Spin representations of symmetric groups and the main theorem}\label{mainthmsec}

\subsection{Representations of symmetric groups and their double covers}\label{snsec}

In this section we summarise the background theory we shall need on representation theory of the symmetric groups and their double covers, sometimes specialising to the case of characteristic $2$, and state our main theorem.

Essential references for the symmetric group are the books of James \cite{jbook} and Mathas \cite{mathbook}; for the double cover $\tsss n$, the book by Hoffman and Humphreys \cite{hh} and (for the case of characteristic $2$) the paper by Bessenrodt and Olsson \cite{bo} are recommended. In contrast to some of these references, we work here mostly with characters rather than modules.

We begin with the ordinary character theory of the symmetric group $\sss n$, which has been well understood for more than a hundred years. For each partition $\la$ of $n$, we let $\spe\la$ denote the \emph{Specht module}, defined over an arbitrary field by James. Over a field of characteristic $0$ the Specht modules are irreducible, and give a complete set of irreducible modules for $\sss n$. We write $\ord\la$ for the character of $\spe\la$.

Now we consider double covers. Let $\tsss n$ denote the double cover of $\sss n$ with generators $s_1,\dots,s_{n-1},z$, subject to relations
\[
z^2=1,\qquad zs_i=s_iz,\qquad s_i^2=1,
\]
\[
s_is_j=zs_js_i\ \ \text{ for }j>i+1,\qquad s_is_js_i=s_js_is_j\ \ \text{ for }j=i+1.
\]
As long as $n\gs4$, $\tsss n$ is a Schur cover of $\sss n$, which means that linear representations of $\tsss n$ are equivalent to projective representations of $\sss n$. If we let $\taaa n$ denote the inverse image of the alternating group $\aaa n$ under the covering map $\tsss n\to\sss n$, then $\taaa n$ is a double cover of $\aaa n$, and is a Schur cover provided $n\gs4$ with $n\neq6,7$.

(In fact for $n\neq6$ there is a different Schur cover of $\sss n$ defined in a similar way, but the two covers are isoclinic, so their representation theories are essentially identical. The double cover $\taaa n$ is the unique Schur cover of $\aaa n$ for $n\neq 6,7$	.)

The ordinary character theory of $\tsss n$ goes back to Schur \cite{schu}. On an irreducible representation (over any field) the central involution $z$ acts as either the identity or minus the identity. In the former case the representation is then just a lift of a representation of $\sss n$; so the characters $\ord\la$ can be extended to irreducible characters of $\tsss n$, which we still write as $\ord\la$. In particular, $\ord{1^n}$ is the one-dimensional sign character, which sends each $s_i$ to $-1$ and $z$ to $1$. Given any class function on $\tsss n$, the \emph{associate} class function is obtained by multiplying by the sign character. 

Representations on which $z$ does not act as the identity are called \emph{spin representations}. The classification of irreducible spin characters can be given as follows. Let us write $\ev\la$ for the number of positive even parts of a partition $\la$. Then for every \trp $\la$ of $n$ with $\ev\la$ even, there is a self-associate irreducible spin character $\lan\la\ran$ of $\tsss n$. For every \trp $\la$ with $\ev\la$ odd, there is a pair of associate irreducible characters $\lan\la\ran_+$, $\lan\la\ran_-$. These characters give a complete set of irreducible spin characters of $\tsss n$. (The characters themselves were found by Schur; explicit constructions of the corresponding representations were given much later, by Nazarov \cite{naz}.)

In this paper we introduce an unusual notation, as follows.

\smallskip
\begin{mdframed}[innerleftmargin=3pt,innerrightmargin=3pt,innertopmargin=3pt,innerbottommargin=3pt,roundcorner=5pt,innermargin=-3pt,outermargin=-3pt]
Suppose $\la\in\cald(n)$. We write
\[
\spn\la=
\begin{cases}
\lan\la\ran&\text{if $\ev\la$ is even}\\
\frac1{\sqrt2}(\lan\la\ran_++\lan\la\ran_-)&\text{if $\ev\la$ is odd}.
\end{cases}
\]
\end{mdframed}

This notation will simplify several formul\ae{} appearing below, including the degree formula for spin characters and the branching rules. Observe that $\spn\la$ is a self-associate class function, and that (with $\chm{\ \ }{\ \ }$ denoting the usual inner product on ordinary characters) $\chm{\spn\la}{\spn\la}=1$.

Now we consider representations in characteristic $2$, again beginning with the symmetric group $\sss n$. If $\la$ is a \trp, then the Specht module $\spe\la$ defined over a field of characteristic $2$ has an irreducible head $\jms\la$, called the \emph{James module}. The modules $\jms\la$ for $\la\in\cald(n)$ give a complete set of irreducible $\sss n$-modules in characteristic $2$. We let $\sid\la$ denote the $2$-modular Brauer character of $\jms\la$.

For the double cover $\tsss n$, observe that there are no irreducible spin characters in characteristic $2$, so a complete set of irreducible $2$-modular Brauer characters is given by the lifts of the characters $\sid\mu$, for $\mu\in\cald(n)$.

For any character $\chi$, we write $\ol\chi$ for its $2$-modular reduction, i.e.\ the $2$-modular Brauer character obtained by taking the values of $\chi$ on elements of odd order.

\subsection{Decomposition numbers}\label{decnosec}

We are interested in computing the decomposition numbers for $\tsss n$; that is, the multiplicities of the characters $\sid\mu$ in the modular reductions of the ordinary irreducible characters. For any ordinary character $\chi$, we write $\cm\chi{\sid\mu}$ for the multiplicity of $\sid\mu$ as an irreducible constituent of $\ol\chi$. In particular, for $\la\in\calp(n)$ we write
\[
D_{\la\mu}=\cm{\ord\la}{\sid\mu}.
\]
The integers $D_{\la\mu}$ are then the entries of the \emph{decomposition matrix} for $\sss n$. For $\la\in\cald(n)$, we write
\[
\dspn_{\la\mu}=\cm{\spn\la}{\sid\mu}.
\]
$\dspn_{\la\mu}$ then belongs to $\bbz$ if $\ev\la$ is even, or to $\sqrt2\bbz$ if $\ev\la$ is odd, so is not quite a decomposition number in the conventional sense. In fact, when $\ev\la$ is odd, the characters $\lan\la\ran_+$ and $\lan\la\ran_-$ differ only on elements of even order, so have the same $2$-modular reduction, with the result that
\[
\cm{\lan\la\ran_+}{\sid\mu}=\cm{\lan\la\ran_-}{\sid\mu}=\tfrac1{\sqrt2}\dspn_{\la\mu}.
\]
So the full decomposition matrix for $\tsss n$ can be recovered immediately  from the numbers $\dspn_{\la\mu}$.

In fact, in characteristic $2$, knowledge of the decomposition matrix of $\sss n$ is enough to determine the decomposition matrix for $\tsss n$. The natural map $\tsss n\to\sss n$ is bijective on (conjugacy classes of) elements of odd order, and via this bijection we can identify $2$-modular Brauer characters for $\sss n$ and $\tsss n$. The $2$-modular reductions $\modr{\ord\xi}$ span the space of $2$-modular Brauer characters, so for any spin character $\spn\la$ we can equate $\mspn\la$ with a linear combination $\sum_\xi a_\xi\modr{\ord\xi}$; in fact, it suffices to use only the characters $\modr{\ord\xi}$ for \trps $\xi$. To find this linear combination, we need to calculate the values of $\spn\la$ and $\ord\xi$ on elements of odd order; there are well-known algorithms for doing this, which are implemented in GAP \cite{gap}. Now given $\mspn\la=\sum_\xi a_\xi\modr{\ord\xi}$, we obtain $\dspn_{\la\mu}=\sum_\xi a_\xi D_{\xi\mu}$ for every $\mu$, and so we have all the decomposition numbers for $\tsss n$.

\subsection{Irreducible and homogeneous characters}\label{rdsec}

Our main interest is in determining which ordinary irreducible characters remain irreducible in characteristic $2$. This problem was solved for $\tsss n$ in \cite{mfspin2}, and in this paper we will solve the same problem for the double cover of the alternating group.

More generally, we will say that an ordinary character $\chi$ is \emph{homogeneous} (in characteristic $2$) if there is only one $\mu$ for which $\cm\chi{\sid\mu}\neq0$. If $\chi$ is an irreducible character of $\sss n$, then $\chi$ being homogeneous in characteristic $2$ is equivalent to it remaining irreducible in characteristic $2$, since (thanks to James \cite{j1}) it is known that $\ol\chi$ has a composition factor that occurs exactly once. But for spin characters of $\tsss n$, the question of homogeneity is more complex than the question of irreducibility. Since we are only concerned with characteristic $2$ in this paper, we will just say ``homogeneous'' to mean homogeneous in characteristic $2$.

In fact, a result of Bessenrodt and Olsson gives us a very good start by identifying a particular composition factor of the $2$-modular reduction of an ordinary irreducible spin character. For a \trp $\la$, recall the partition $\la\dblreg$ defined in \cref{regdoubsec}.

\begin{thmc}{bo}{Theorem 5.2}\label{spinreg}
Suppose $\la\in\cald$. Then $\decs\la{(\la\dblreg)}=2^{\ev\la/2}$.
\end{thmc}

This is very useful, because it means that if $\spn\la$ is homogeneous, then $\ol{\spn\la}$ must equal $2^{\ev\la/2}\sid{\la\dblreg}$. Hence an irreducible character $\lan\la\ran$ or $\lan\la\ran_{\pm}$ is irreducible in characteristic $2$ \iff $\spn\la$ is homogeneous and $\ev\la$ equals $0$ or $1$.

\subsection{The double cover of the alternating group and the main theorem}\label{ansec}

The main theorem in this paper concerns the double cover $\taaa n$ of the alternating group. As with $\tsss n$, the ordinary irreducible characters of $\taaa n$ come in two types: those lifted from characters of $\aaa n$, and spin characters. For characters lifted from $\aaa n$, the question of irreducibility in characteristic $2$ is exactly the same as for $\aaa n$, and this question was answered in \cite{mfonaltred}. So in this paper we need only address spin characters. The ordinary irreducible spin characters of $\aaa n$ are labelled by \trps $\la$ of $n$, but here the situation is opposite to that in $\tsss n$: if $\ev\la$ is even, there is a pair of associate spin characters $\lan\la\ran_+$, $\lan\la\ran_-$, while if $\ev\la$ is odd there is a self-associate spin character $\lan\la\ran$. These characters are obtained by restricting the corresponding characters of $\tsss n$, and give a complete set of ordinary irreducible spin characters of $\aaa n$.

In \cite{mfspin2} our main question for $\taaa n$ was reduced to a problem about $\tsss n$, as follows.

\begin{propnc}{mfspin2}{Section 7.1}\label{sssaaareduction}
Suppose $\la\in\cald(n)$. Then the spin character $\lan\la\ran$ or $\lan\la\ran_\pm$ for $\taaa n$ remains irreducible in characteristic $2$ \iff $\spn\la$ is homogeneous, $\ev\la\ls2$ and $\la$ does not have the form $(4b)$ for $b\gs1$.
\end{propnc}

We are now almost ready to state our main theorem. First we need to recall the definition of a $2$-Carter partition. Say that $\la\in\calp$ is \emph{$2$-Carter} if for every $r\gs1$, $\la_r-\la_{r+1}+1$ is divisible by a power of $2$ greater than~$\la_{r+1}-\la_{r+2}$.

An equivalent definition of $2$-Carter partitions may be given in terms of hook lengths: if $(r,c)$ is a node of $\la$, the \emph{$(r,c)$-hook length} $\hkl rc\la$ is the integer $\la_r+\la'_c-r-c+1$. It is proved in \cite[Lemma 3.14]{jsome} that $\la$ is $2$-Carter \iff $\nu_2(\hkl rc\la)=\nu_2(\hkl sc\la)$ for every pair of nodes $(r,c)$ and $(s,c)$ in the same column of $[\la]$.

From \cref{sssaaareduction} the character(s) labelled by $\la\in\cald(n)$ can be irreducible in characteristic $2$ only if $\ev\la\ls2$. For the cases where $\ev\la\ls1$, we can read the classification directly from \cite[Theorem 3.3]{mfspin2}, using \cref{sssaaareduction}.

\begin{thm}\label{mainthm01}
Suppose $\la\in\cald(n)$.
\begin{enumerate}
\item
If $\ev\la=0$, then the spin characters $\lan\la\ran_+$ and $\lan\la\ran_-$ of $\taaa n$ remain irreducible in characteristic $2$ \iff $\la$ has the form $\tau+4\alpha$, where $\tau$ is a \fbc and $\alpha$ is a $2$-Carter partition with $\len\tau\gs\len\alpha$.
\item
If $\ev\la=1$, then the spin character $\lan\la\ran$ of $\taaa n$ remains irreducible in characteristic $2$ \iff one of the following occurs:
\begin{enumerate}
\item
$\la$ has the form $\tau+4\alpha\sqcup(2)$, where $\tau$ is a \fbc and $\alpha$ is a $2$-Carter partition with $\tau_{\len\alpha}\gs3$;
\item
$\la$ equals $(4b-2)$ or $(4b-2,1)$ for some $b\gs2$;
\item
$\la=(3,2,1)$.
\end{enumerate}
\end{enumerate}
\end{thm}

So our task in the present paper is to address the case where $\ev\la=2$. Our main result is the following.

\smallskip
\begin{mdframed}[innerleftmargin=3pt,innerrightmargin=3pt,innertopmargin=3pt,innerbottommargin=3pt,roundcorner=5pt,innermargin=-3pt,outermargin=-3pt]
\vspace*{-\topsep}
\begin{thm}\label{mainthm}
Suppose $\la\in\cald(n)$ with $\ev\la=2$. Then the irreducible spin characters $\lan\la\ran_+$ and $\lan\la\ran_-$ of $\taaa n$ remain irreducible in characteristic $2$ \iff one of the following occurs:
\begin{enumerate}
\item
$\la$ has the form $\tau+4\alpha\sqcup(4,2)$ where $\tau$ is a \fbc, $\alpha$ is a $2$-Carter partition and $\len\tau>\len\alpha$;
\item
$\la$ is one of $(4,3,2,1)$, $(5,4,3,2)$, $(5,4,3,2,1)$, $(7,4,3,2,1)$.
\end{enumerate}
\end{thm}
\end{mdframed}

For the rest of the paper we will work with $\tsss n$; in view of \cref{sssaaareduction}, our task is to show that for $\la\in\cald(n)$ with $\ev\la=2$, the character $\spn\la$ for $\tsss n$ is homogeneous \iff one of the conditions in \cref{mainthm} occurs.

\medskip
It is an interesting question to ask when $\spn\la$ is homogeneous in general. In \cref{rouqsec} we answer this for the case where $\la$ is a \emph{separated} partition, which appears to be the generic case. A full solution involves understanding how the sequence of ``exceptional'' partitions from \cref{mainthm01,mainthm}, namely
\[
(3,2,1),\ (4,3,2,1),\ (5,4,3,2),\ (5,4,3,2,1),\ (7,4,3,2,1),\ \dots
\]
continues. It appears from calculations that the partitions in this sequence are close to being $2$-cores. In particular, we conjecture that if $\la$ is a $2$-core, then $\spn\la$ is homogeneous. But the more general pattern remains mysterious.

\medskip
For the remainder of \cref{mainthmsec} we recall some more of the background we shall need on the representation theory of $\sss n$ and $\tsss n$, and prove some preliminary results.

\subsection{The Schur algebra and the $(-1)$-Schur algebra}\label{schursec}

In this paper we shall need to consider the Schur algebra introduced by Green \cite{gbook}, and its quantum analogue introduced by Dipper and James \cite{dj}. We let $\schur n$ denote the classical Schur algebra defined over an infinite field $\bbf$ of characteristic $2$. We also let $\qschur n$ denote the $q$-Schur algebra over $\bbc$ with quantum parameter $q=-1$. These algebras are denoted $\cals_\bbf(1,n)$ and $\cals_\bbc(-1,n)$ in the notation of Dipper and James.

For any partition $\la\in\calp(n)$, there is a \emph{Weyl module} $\weyl\la$ for $\schur n$, with a simple head $\schir\la$. The modules $\schir\la$ give a complete set of irreducible modules for $\schur n$. Note that we use the convention from James's paper \cite{j10}, in which the Specht module $\spe\la$ is the image of $\weyl\la$ (not $\weyl{\la'}$) under the Schur functor.

Given $\la,\mu\in\calp(n)$, we write $D_{\la\mu}$ for the composition multiplicity $[\weyl\la:\schir\mu]$ for the Schur algebra $\schur n$. This is consistent with the notation  introduced in \cref{decnosec}, since $[\spe\la:\sid\mu]=[\weyl\la:\schir\mu]$ when $\mu$ is \treg.

Correspondingly, there are Weyl modules $\weyl\la$ with simple heads $\schir\la$ for the $(-1)$-Schur algebra $\qschur n$. In this case, we write $\odc_{\la\mu}$ for the decomposition number $[\weyl\la:\schir\mu]$. These decomposition numbers are better understood than the decomposition numbers for $\schur n$ in that there is an algorithm to compute them. However, the two sets of decomposition numbers are closely related. For a fixed $n$, let $D$ denote the matrix with entries $D_{\la\mu}$ for $\la,\mu\in\calp(n)$, and define $\odc$ similarly.

\begin{propnc}{mathbook}{Theorem 6.35}\label{adjmat}
For a given $n$, there is a square matrix $A$ with non-negative integer entries such that $D=\odc A$.
\end{propnc}

The matrix $A$ is called the \emph{adjustment matrix} for $\schur n$. The decomposition numbers for both $\schur n$ and $\qschur n$ satisfy the following unitriangularity property, which comes from the fact that $q$-Schur algebras are quasi-hereditary.

\begin{lemma}\label{dndom}
Suppose $\la,\mu\in\calp(n)$. Then $D_{\la\la}=\odc_{\la\la}=1$, while $D_{\la\mu}=\odc_{\la\mu}=0$ unless $\mu\dom\la$.
\end{lemma}

As a consequence, the adjustment matrix is also unitriangular, i.e.\ $A_{\la\la}=1$, while $A_{\la\mu}=0$ when $\mu\ndom\la$.

A key component in the proof of our main theorem will be the classification of irreducible Weyl modules.

\begin{thmc}{jm1}{Theorem 4.5}\label{irredweyl}
Suppose $\la\in\calp(n)$. Then the Weyl module $\weyl\la$ for $\schur n$ is irreducible \iff $\la$ is $2$-Carter.
\end{thmc}

(Note that \cite{jm1} uses the opposite convention for labelling Weyl modules to ours, so partitions need to be conjugated when reading \cref{irredweyl} from that reference.)

\subsection{Blocks}\label{blocksec}

In the modular representation theory of any group, it is useful to sort characters into blocks. Here we summarise the combinatorics underlying the $2$-blocks of $\sss n$ and $\tsss n$. Since we are only concerned with characteristic $2$ in this paper, we will say simply ``block'' to mean ``$2$-block''. The material here is taken directly from \cite{mfspin2}.

Recall from \cref{corequotsec} that the $2$-core of a partition $\la$ is the partition obtained by repeatedly removing rim $2$-hooks. The number of rim $2$-hooks removed we call the \emph{$2$-weight} of $\la$. The following result is a special case of the Brauer--Robinson Theorem \cite{brauer,gdbr}.

\begin{thm}\label{brarob}
Suppose $\la,\mu\in\calp(n)$. Then $\ord\la$ and $\ord\mu$ lie in the same $2$-block of $\sss n$ \iff $\la$ and $\mu$ have the same $2$-core.
\end{thm}

The distribution of irreducible Brauer characters into blocks is given by the same rule: since a Brauer character $\sid\mu$ occurs as a composition factor of $\ord\mu$, it lies in the same block as $\ord\mu$.

\cref{brarob} yields a convenient way of labelling the $2$-blocks of $\sss n$: we can define the core of a block to be the common $2$-core of the partitions labelling the irreducible characters in that block. Two partitions with the same size and the same $2$-core obviously have the same $2$-weight, so if $n$ is unspecified we can label a block by its core and weight, i.e.\ the common $2$-core and $2$-weight of the partitions labelling characters in that block.

It is sometimes also convenient to label blocks in terms of residues: recall from \cref{corecont} that two partitions of the same size have the same $2$-core \iff they have the same content. So we can define the content of a block to be the content of any partition labelling an irreducible character in that block.

A similar statement \cite[Theorem 4.24]{jm1} applies for the blocks of the Schur algebra $\schur n$ and the $(-1)$-Schur algebra $\qschur n$: two Weyl modules $\weyl\la$ and $\weyl\mu$, or two simple modules $\schir\la$ and $\schir\mu$, lie in the same block (for either algebra) \iff $\la$ and $\mu$ have the same $2$-core.

\smallskip Now we come to the $2$-blocks of $\tsss n$; these were determined by Bessenrodt and Olsson \cite[Theorem 4.1]{bo}. For the lifts of the ordinary characters of $\sss n$, it is still the case that $\ord\la$ and $\ord\mu$ lie in the same $2$-block \iff $\la$ and $\mu$ have the same $2$-core; the same statement holds also for the lifts of the Brauer characters $\sid\mu$. Since these are the only irreducible Brauer characters of $\tsss n$, it follows that there are no other $2$-blocks. So the $2$-blocks of $\tsss n$ can also be labelled by their $2$-cores, and to complete the picture of $2$-blocks it suffices to say which blocks the spin characters $\spn\la$ lie in. (This question makes sense even with our unusual definition of $\spn\la$, because when $\ev\la$ is odd the characters $\lan\la\ran_+$ and $\lan\la\ran_-$ lie in the same block.) In fact this can be inferred directly from \cref{spinreg}: $\spn\la$ lies in the block whose core is the $2$-core of $\la\dblreg$.

There is a more direct way to see whether two spin characters lie in the same block: if $\la,\mu\in\cald(n)$, then by \cref{regdoub} $\spn\la$ and $\spn\mu$ lie in the same block of $\tsss n$ \iff $\la$ and $\mu$ have the same \fbc. So we can define the \emph{\bc} of a block of $\tsss n$ to be the common \fbc of the partitions labelling spin characters in that block.

We can also express the block classification for spin characters in terms of residues. We define the \emph{\spr} of a node $(r,c)$ to be the residue of $\lfloor c/2\rfloor$ modulo $2$. Then we can define the spin-content of a \trp to be the multiset of the \sprs of its nodes. It is an easy exercise to check that the spin-content of $\la$ equals the content of $\la\dblreg$, so $\spn\la$ lies in the block whose content is the spin-content of $\la$.

For example, take $\la=(6,2,1)$, for which $\la\dbl=(3^2,1^3)$ and $\la\dblreg=(5,3,1)$. The character $\spn\la$ lies in the block with content $\{0^5,1^4\}$, as we see from the following diagram for these three partitions. In the diagram for $\la$, nodes are labelled with their spin-residues, while in the other two diagrams nodes are labelled with their residues.
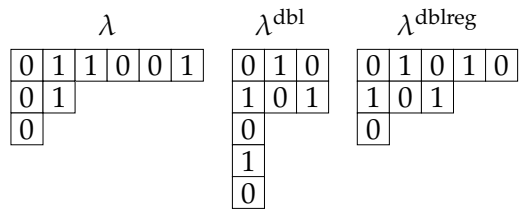
\[
\begin{array}{ccc}
\la&\la\dbl&\la\dblreg
\\
\young(011001,01,0)&\young(010,101,0,1,0)&\young(01010,101,0)
\end{array}
\]

\subsection{The inverse of the decomposition matrix of the Schur algebra}

A very useful result due to James is a reduction theorem for the entries of the inverse of the decomposition matrix for a $q$-Schur algebra; this derives from Steinberg's tensor product theorem. The cases we need are as follows.

\begin{propn}\label{steinberggeneral}
Suppose $\alpha,\mu\in\calp$ and $\beta\in\cald$ with $2|\alpha|+|\beta|=|\mu|$. Then
\begin{align*}
\odc\v_{(\dup\alpha\sqcup\beta)\mu}&=\sum_{\la,\xi\in\calp}\odc\v_{\beta\la}\epstar\xi\lrc\alpha{\quo\xi0}{\quo\xi1}\lrc\mu\la\xi
\\
\intertext{and}
D\v_{(\dup\alpha\sqcup\beta)\mu}&=\sum_{\la,\xi,\gamma\in\calp}D\v_{\beta\la}D\v_{\alpha\gamma}\epstar\xi\lrc\gamma{\quo\xi0}{\quo\xi1}\lrc\mu\la\xi
\end{align*}
\end{propn}

\begin{pf}
This is proved in \cite[Corollary 6.9]{j10}, but quite a bit of translation of notation is needed, so we explain this here. If we take the first formula in \cite{j10} with $e=p=2$, then the terms $c$ and $c^\ast$ are both just our $D\v$. Replacing the variables $\rho,\alpha,\beta,\sigma$ with $\alpha,\beta,\la,\gamma$ in that formula, we obtain
\[
D\v_{(\dup\alpha\sqcup\beta)\mu}=(-1)^{\card\alpha}\sum_{\la,\gamma\in\calp}D\v_{\beta\la}D\v_{\alpha\gamma}\upsilon(\la,\gamma,\mu).
\]
The term $\upsilon(\la,\gamma,\mu)$ is defined in \cite[Definitions 2.13 and 2.20]{j10} as
\[
\sum_\xi\epsilon(\xi)\lrc\gamma{\quo\xi0}{\quo\xi1}\lrc\mu\la\xi,
\]
summing over $\xi\in\calp(2\card\gamma)$ with empty $2$-core. The term $\epsilon(\xi)$ is defined like our $\epstar\xi$ but with horizontal and vertical rim-hooks interchanged, so that $\epsilon(\xi)=(-1)^{\card\gamma}\epstar\xi$. Substituting this above, we obtain our second formula.

For our first formula, we need to take $e=2$ and $p$ very large, and use the second formula in \cite{j10}. (It is known that the decomposition matrix for the $(-1)$-Schur algebra over $\bbc$ is the same as over a field of very large finite characteristic.) Now $c$ is our $\odc\v$, so the formula in \cite{j10} gives our first formula.
\end{pf}

\begin{cory}\label{steincory}
Suppose $\alpha,\mu\in\calp$ and $\beta$ is a $2$-core with $2|\alpha|+|\beta|=|\mu|$. Then
\begin{align*}
\odc\v_{(\dup\alpha\sqcup\beta)\mu}&=\sum_{\xi\in\calp}\epstar\xi\lrc\alpha{\quo\xi0}{\quo\xi1}\lrc\mu\beta\xi,
\\
D\v_{(\dup\alpha\sqcup\beta)\mu}&=\sum_{\gamma\in\calp}D\v_{\alpha\gamma}\odc\v_{(\dup\gamma\sqcup\beta)\mu},
\\
A\v_{(\dup\alpha\sqcup\beta)\mu}&=
\begin{cases}
D\v_{\alpha\gamma}&\text{if $\mu=\dup\gamma\sqcup\beta$ for some $\gamma\in\calp$}\\
0&\text{otherwise},
\end{cases}
\\
A_{(\dup\alpha\sqcup\beta)\mu}&=
\begin{cases}
D_{\alpha\gamma}&\text{if $\mu=\dup\gamma\sqcup\beta$ for some $\gamma\in\calp$}\\
0&\text{otherwise}.
\end{cases}
\end{align*}
\end{cory}

\begin{pf}
This \lcnamecref{steincory} was proved in the special cases $\beta=\varnothing,(1)$ in \cite{mfspin2}. The fact that $\beta$ is a $2$-core means that for both $\schur{\card\beta}$ and $\qschur{\card\beta}$ the Weyl module $\weyl\beta$ lies in a block by itself, so that
\[
\odc\v_{\beta\la}=D\v_{\beta\la}=\delta_{\beta\la}.
\]
Substituting this into the two formul\ae{} in \cref{steinberggeneral} gives the first two statements. The last two statements are proved exactly as in \cite[Corollary 2.5]{mfspin2}.
\end{pf}

\subsection{Induction and restriction}\label{brnchsec}

It will be very helpful for us to use induction and restriction between different symmetric groups and their double covers. We recall the essential background we shall need; much of this is taken from \cite{mfspin2}, where further references can be found.

Given a character $\chi$ of $\tsss n$, we write $\chi\res_{\tsss{n-1}}$ for its restriction to $\tsss{n-1}$, and $\chi\ind^{\tsss{n+1}}$ for the corresponding induced character for $\tsss{n+1}$. In fact, we use refinements of these operations, introduced in the symmetric group case by Robinson. Suppose $\chi$ is a character of $\tsss n$, lying in the block $B$ with content $\muset ab$. Then we write $\ee0\chi$ for the component of $\chi\res_{\tsss{n-1}}$ lying in the block with content $\muset{a-1}b$ if there is such a block, and set $\ee0\chi=0$ otherwise. Similarly, we write $\ee1\chi$ for the component of $\chi\res_{\tsss{n-1}}$ lying in the block with content $\muset a{b-1}$ if there is such a block, and set $\ee1\chi=0$ otherwise. We extend the functions $\ee0,\ee1$ linearly. These functions $\ee0,\ee1$ can be applied either to ordinary characters or to $2$-modular Brauer characters, and for any character $\chi$ we have $\chi\res_{\tsss{n-1}}=\ee0\chi+\ee1\chi$ (this follows from the classical branching rule for ordinary irreducible representations of $\sss n$, together with the block classification). Defining these functions for any $n$, we can consider powers $\ee i^a$ for $a\gs0$. In fact, it will be useful to define \emph{divided powers} $\eed ia=\ee i^a/a!$. Given a non-zero character $\chi$ and a residue $i$, we define $\epsilon_i\chi$ to be the largest $a\gs0$ for which $\eed ia\neq0$, and we write $\emx i\chi=\eed i{\epsilon_i\chi}\chi$.

A similar situation applies for induction of characters to $\tsss{n+1}$: we can write $\chi\ind^{\tsss{n+1}}=\ff0\chi+\ff1\chi$, where $\ff0$ and $\ff1$ are functions defined using the block classification in a similar way to $\ee0$ and $\ee1$. We define divided powers $\ffd ia$, and for a non-zero character $\chi$ we define $\phi_i\chi$ to be the largest $a$ for which $\ffd ia\chi$ is non-zero, and write $\fmx i\chi=\ffd i{\phi_i\chi}\chi$.

We can describe the effect of these functions on the irreducible spin characters as follows. Given \trps $\la$ and $\mu$, a residue $i\in\{0,1\}$ and an integer $a\gs0$, we write $\la\ads ia\mu$ if $\mu$ can be obtained from $\la$ by adding $a$ nodes of \spr $i$. Then we have the following, which extends the spin branching rule of Dehuai and Wybourne \cite{dw}.

\begin{propn}\label{spinbranchpower}
Suppose $\la\in\cald(n)$, $\mu\in\cald(n-a)$ and $i\in\{0,1\}$. Then $\chm{\eed ia\spn\la}{\spn\mu}>0$ \iff $\mu\ads ia\la$. If this is the case, let $c$ be the number of pairs of nodes of $\la\setminus\mu$ lying in adjacent columns. Then
\[
\chm{\eed ia\spn\la}{\spn\mu}=2^{(a-\len\la+\len\mu)/2-c}.
\]
\end{propn}

\begin{pf}
This is essentially proved in \cite[Proposition 2.17]{mfspin2}, though the result there is more complicated than the formula here, with four cases depending on the parities of $\ev\la$ and $\ev\mu$. The present version is simpler because of our definition of the notation $\spn\la$.
\end{pf}

An analogous result holds for $\ffd ia\spn\la$, and we refer to these results together as the \emph{spin branching rules}.

From the spin branching rule we can immediately determine the effect of the functions $\emx i$ and $\fmx i$ on spin characters. Let $\la$ be a \trp, and let $\la\dn i$ be the smallest \trp which can be obtained by removing nodes of \spr $i$ from $\la$. We refer to these nodes as the \emph{$i$-\sprms} of $\la$. From the spin branching rule, $\emx i\spn\la$ is a non-zero scalar multiple of $\spn{\la\dn i}$. Similarly, we define $\la\up i$ to be the largest \trp which can be obtained by adding nodes of \spr $i$ to $\la$. These nodes are called the \emph{$i$-\spams} of $\la$, and $\fmx i\spn\la$ is a non-zero multiple of $\spn{\la\up i}$.

We define a \sprm of $\la$ to be a node which is either a $0$- or a $1$-\sprm of $\la$; in other words, a node of $\la$ that can be removed (possibly together with some other nodes of the same \spr) to leave a smaller \trp. Note that a \sprm is not necessarily a removable node as defined in \cref{compnpartnsec}; for example, the node $(2,2)$ is a \sprm of the partition $(4,3)$. Also, a node which is removable in the conventional sense might not be spin-removable if removing it would cause the partition no longer to be \treg. We define \spams similarly to \sprms.

\medskip
We will also need to recall Kleshchev's modular branching rules, which partially describe the effect of $\ee i$ and $\ff i$ on irreducible $2$-modular Brauer characters. In fact, we will just need to know the effect of $\emx i$ and $\fmx i$. To set this up, we need to consider sign sequences. Suppose $s=s_1\dots s_r$ is a finite sequence consisting of signs $+$ and $-$. The \emph{reduction} of $s$ is the subsequence obtained by repeatedly deleting successive pairs $+-$.

Now suppose $\mu$ is a \trp and $i\in\{0,1\}$. The \emph{$i$-signature} of $\mu$ is defined to be the sign sequence obtained by working from top to bottom of the Young diagram of $\mu$, writing a $+$ for each addable $i$-node and a $-$ for each removable $i$-node. The reduction of this sign sequence is called the \emph{reduced $i$-signature} of $\mu$. The removable nodes corresponding to minus signs in the reduced $i$-signature are called the \emph{normal} $i$-nodes of $\mu$, and the addable nodes corresponding to the plus signs in the reduced $i$-signature are called the \emph{conormal} $i$-nodes.

With these definitions, we have the following. This (and many more results) can be found in the survey \cite{bk} (in particular, see the discussion following Lemma 2.12).

\begin{thm}\label{modbranch}
Suppose $\chi$ is an irreducible $2$-modular Brauer character of $\tsss n$, and $i\in\{0,1\}$. Then $\emx i\chi$ and $\fmx i\chi$ are irreducible Brauer characters. Specifically, write $\chi=\sid\mu$ for $\mu\in\cald$, and let $\mu^-$ be the partition obtained by removing all the normal $i$-nodes from $\mu$, and $\mu^+$ the partition obtained by adding all the conormal $i$-nodes to $\mu$. Then $\mu^-,\mu^+$ are \treg, and
\[
\emx i\sid\mu=\sid{\mu^-},\qquad\fmx i\sid\mu=\sid{\mu^+}.
\]
\end{thm}

Now take a \trp $\la$, and recall that $\sid{\la\dblreg}$ is an irreducible constituent of the Brauer character $\spn\la$. Since restriction is an exact functor, this immediately gives $\epsilon_i\sid{\la\dblreg}\ls\epsilon_i\spn\la$, and we write $\noregdown i\la$ if $\epsilon_i\sid{\la\dblreg}<\epsilon_i\spn\la$. It will be very useful to be able to determine from $\la$ (without calculating $\la\dblreg$) exactly when $\noregdown i\la$. In fact a complete answer to this question is quite complicated, so we just give a simple sufficient condition for $\noregdown i\la$.

First we need some results on counting nodes and addable and removable nodes. Recall that for $l\gs0$ the $l$th ladder in $\bbn^2$ is the set of nodes $(r,c)$ with $r+c=l+2$. Given a partition $\mu$, we write
\begin{align*}
\lad l\mu\ &\text{for the number of nodes of $\mu$ in ladder $l$,}\\
\alad l\mu\ &\text{for the number of addable nodes of $\mu$ in ladder $l$,}\\
\rlad l\mu\ &\text{for the number of removable nodes of $\mu$ in ladder $l$,}
\end{align*}
setting all of these numbers to be zero when $l<0$.

These numbers are related to each other by the following result. (Here and throughout, we use the Kronecker delta.)

\begin{lemma}\label{ladaddrem}
If $\mu\in\calp$ and $l\gs0$, then
\[
\alad l\mu-\rlad{l-2}\mu=\delta_{l0}-\lad l\mu+2\lad{l-1}\mu-\lad{l-2}\mu.
\]
\end{lemma}

\begin{pf}
A version of this lemma for arbitrary characteristic is given in \cite[Lemma 4.8]{mfaltred}, though it is only proved there in odd characteristic (which is the only case needed in that paper). We prove it by induction on $|\mu|$, with the case $\mu=\varnothing$ being trivial. Assuming $\mu\neq\varnothing$, let $(r,c)$ be a removable node of $\mu$, and let $m=r+c-2$. Let $\xi$ be the partition obtained from $\mu$ by removing the node $(r,c)$. Then for any $k$ (writing $\indi S$ for the indicator function of the truth of a statement $S$)
\begin{alignat*}2
\lad k\mu&=\lad k\xi&&+\indi{k=m}\\
\alad k\mu&=\alad k\xi&&-\indi{k=m}\\
&&&+\indi{k=m+1\text{ and either }r=1\text{ or }(r-1,c+1)\in\xi}\\
&&&+\indi{k=m+1\text{ and either }c=1\text{ or }(r+1,c-1)\in\xi}\\
\rlad k\mu&=\rlad k\xi&&+\indi{k=m}\\
&&&-\indi{k=m-1\text{, }r>1\text{ and }(r-1,c+1)\notin\xi}\\
&&&-\indi{k=m-1\text{, }c>1\text{ and }(r+1,c-1)\notin\xi}.
\end{alignat*}
Hence the result of the \lcnamecref{ladaddrem} holds for $\mu$ \iff it holds for $\xi$.
\end{pf}

We now give a spin analogue of \cref{ladaddrem}. For $l\gs0$, define the $l$th \emph{slope} in $\bbn^2$ to be the set of nodes $(r,c)$ for which $2r+\inp{c/2}=l+2$. We say that the $m$th slope is \emph{longer} than the $l$th slope if $m>l$. Given a \trp $\la$, we write
\begin{align*}
\slp l\la\ &\text{for the number of nodes of $\la$ in slope $l$,}\\
\aslp l\la\ &\text{for the number of \spams of $\la$ in slope $l$,}\\
\rslp l\la\ &\text{for the number of \sprms of $\la$ in slope $l$,}
\end{align*}
again setting all these numbers to be zero when $l<0$. We also need to define $\zslp l\la$ to be the number of nodes $(r,c)$ in slope $l$ such that $c$ is even, $r\gs2$, and $(\la_{r-1},\la_r,\la_{r+1})=(c+1,c,c-1)$.

\begin{eg}
Take $\la=(15,9,8,7,1)$ and $l=8$. The Young diagram of $\la$ with nodes in the $8$th slope marked is as follows.
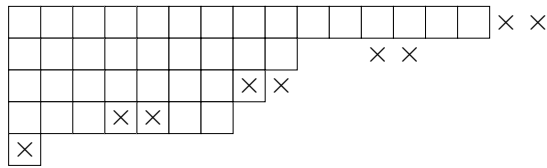
\[
\gyoung(;;;;;;;;;;;;;;;:\tim:\tim,;;;;;;;;;:::\tim:\tim,;;;;;;;;\tim:\tim,;;;;\tim;\tim;;,;\tim)
\]
We see that $\slp6\la=4$, $\aslp6\la=2$ (note that both of the nodes $(1,16)$ and $(1,17)$ are spin-addable, but $(3,9)$ is not) and $\rslp6\la=1$. Also $\zslp6\la=1$, because of the node $(3,8)$.
\end{eg}

Ladders and slopes are related by the following \lcnamecref{ladslope}, which is stated implicitly in \cite{bo}.

\begin{lemmac}{mfspin2}{Lemma 2.11}\label{ladslope}
Suppose $\la$ is a \trp. Then for every $l\gs0$,
\[
\slp l\la=\lad l{\la\dblreg}.
\]
\end{lemmac}

Now we give the analogue of \cref{ladaddrem} for slopes.

\begin{lemma}\label{slpaddrem}
If $\la\in\cald$ and $l\gs0$, then
\[
\aslp l\la-\rslp{l-2}\la=\delta_{l0}-\slp l\la+2\slp{l-1}\la-\slp{l-2}\la+\zslp l\la-\zslp{l-2}\la.
\]
\end{lemma}

\begin{pf}
Suppose $(r,c)$ is a node in slope $l$, with $c$ odd. Assuming $r\gs3$ and $c\gs3$, consider the set of eight nodes
\[
(r-2,c+1),(r-1,c-1),(r-1,c),(r-1,c+1),(r,c-2),(r,c-1),(r,c),(r+1,c-2)
\]
lying in slopes $l-3,l-2,l-2,l-1,l-1,l,l,l+1$ respectively. The configuration of these nodes (labelled with the slopes containing them) is as follows.
{\small
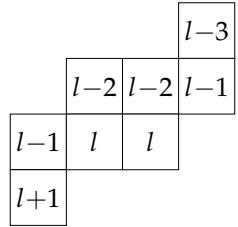
\[
\gyoungs(1.75,:::;<l{-}3>,:;<l{-}2>;<l{-}2>;<l{-}1>,;<l{-}1>;l;l,;<l{+}1>)
\]
}
By considering the possibilities for which of these nodes are nodes of $\la$, we can check that the formula in the \lcnamecref{slpaddrem} is true when restricted to these eight nodes. The possible cases are given in the following table, where for each case we shade the nodes contained in $\la$ and give the contributions to the formula coming from these eight nodes.
{\Yvcentermath1\Yboxdim{8pt}
\begin{longtable}{cccccccc}\hline
&$\slp{l-2}\la$&$\slp{l-1}\la$&$\slp l\la$&$\rslp{l-2}\la$&$\aslp l\la$&$\zslp{l-2}\la$&$\zslp l\la$
\\\hline\\[-9pt]\endhead
\gyoung(:::;,:;;;,;;;,;)
&0&0&0&0&0&0&0
\\
\gyoung(:::!\gry;,:!\wht;;;,;;;,;)
&0&0&0&0&0&0&0
\\
\gyoung(:::;,:!\gry;!\wht;;,;;;,;)
&1&0&0&1&0&0&0
\\
\gyoung(:::!\gry;,:;!\wht;;,;;;,;)
&1&0&0&1&0&0&0
\\
\gyoung(:::;,:!\gry;!\wht;;,!\gry;!\wht;;,;)
&1&1&0&0&0&1&0
\\
\gyoung(:::!\gry;,:;;!\wht;,;;;,;)
&2&0&0&2&0&0&0
\\
\gyoung(:::!\gry;,:;!\wht;;,!\gry;!\wht;;,;)
&1&1&0&0&1&0&0
\\
\gyoung(:::!\gry;,:;;;!\wht,;;;,;)
&2&1&0&0&0&0&0
\\
\gyoung(:::!\gry;,:;;!\wht;,!\gry;!\wht;;,;)
&2&1&0&1&1&0&0
\\
\gyoung(:::!\gry;,:;;;,;!\wht;;,;)
&2&2&0&0&2&0&0
\\
\gyoung(:::!\gry;,:;;!\wht;,!\gry;;!\wht;,;)
&2&1&1&1&0&0&0
\\
\gyoung(:::!\gry;,:;;;,;;!\wht;,;)
&2&2&1&0&1&0&0
\\
\gyoung(:::!\gry;,:;;!\wht;,!\gry;;!\wht;,!\gry;!\wht)
&2&1&1&0&0&0&1
\\
\gyoung(:::!\gry;,:;;;,;;!\wht;,!\gry;!\wht)
&2&2&1&0&1&0&0
\\
\gyoung(:::!\gry;,:;;;,;;;,!\wht;)
&2&2&2&0&0&0&0
\\
\gyoung(:::!\gry;,:;;;,;;;,;!\wht)
&2&2&2&0&0&0&0
\end{longtable}
}
A similar statement holds if $r\ls2$ or $c=1$ (where there are fewer nodes to consider). Summing over all pairs $(r,c)$, we account for the nodes in slopes $l-2,l-1,l$ once each, and obtain the desired result.
\end{pf}

Before giving our main result, we need a simple lemma about reduction of sign sequences. We leave the proof as an exercise.

\begin{lemma}\label{signseq}
Suppose $s=s_1\dots s_r$ is a sign sequence. Let $m$ denote the total number of minus signs in $s$. Suppose that for some $t$ we can find distinct integers $a_1,\dots,a_t,b_1,\dots,b_t$ such that $s_{a_i}=+$, $s_{b_i}=-$ and $a_i<b_i$ for each $i$. Then the total number of minus signs in the reduction of $s$ is at most $m-t$.
\end{lemma}

Finally we can give our sufficient condition to have $\noregdown i\la$ when $\la\in\cald$. This is analogous to \cite[Proposition 4.9]{mfaltred} for linear representations of symmetric groups, though weaker, in that the condition we give is sufficient but not necessary for $\noregdown i\la$.

\begin{propn}\label{noregcond}
Suppose $\la$ is a \trp and $i\in\{0,1\}$, and that $\la$ has both an $i$-\sprm and an $i$-\spam, with the $i$-\spam in a longer slope than the $i$-\sprm. Then $\noregdown i\la$.
\end{propn}

\begin{pf}
For this proof we use the following notation: given an integer $u$, define $\sse u$ to be a sequence of plus signs of length $u$ if $u\gs0$, or a sequence of minus signs of length $-u$ if $u<0$. Observe that if $u,v\in\bbz$ with either $u\gs0$ or $v\ls0$, then the reduction of the concatenation $\sse u\sse v$ is $\sse{u+v}$.

We write $\mu=\la\dblreg$. From the remarks following the proof of \cref{spinbranchpower}, we have
\[
\epsilon_i\spn\la=\sum_{l\equiv i\ppmod2}\rslp l\la.
\]
We need to show that $\epsilon_i\sid\mu$ is less than this, and we consider the $i$-signature of $\mu$. Observe that as we read the addable and removable $i$-nodes of a \trp from top to bottom, nodes in higher-numbered ladders occur before nodes in lower-numbered ladders, and within any ladder the removable nodes come before the addable nodes. So the $i$-signature of $\mu$ is
\[
\dots\sse{\alad{i+6}\mu}\sse{-\rlad{i+4}\mu}\sse{\alad{i+4}\mu}\sse{-\rlad{i+2}\mu}\sse{\alad{i+2}\mu}\sse{-\rlad i\mu}\sse{\alad i\mu}.
\]
By the observation above, the reduction of this sequence is the same as the reduction of
\[
\dots\sse{\alad{i+6}\mu-\rlad{i+4}\mu}\sse{\alad{i+4}\mu-\rlad{i+2}\mu}\sse{\alad{i+2}\mu-\rlad i\mu}\sse{\alad i\mu},
\]
and by Lemmas \ref{ladaddrem}, \ref{ladslope} and \ref{slpaddrem} this equals
\begin{align*}
\dots{}&\sse{\aslp{i+6}\la-\rslp{i+4}\la-\zslp{i+6}\la+\zslp{i+4}\la}
\\
&\sse{\aslp{i+4}\la-\rslp{i+2}\la-\zslp{i+4}\la+\zslp{i+2}\la}
\\
&\sse{\aslp{i+2}\la-\rslp{i}\la-\zslp{i+2}\la+\zslp{i}\la}
\\
&\sse{\aslp{i}\la-\zslp{i}\la}.
\end{align*}
Again using the observation above, the reduction of this sequence is the same as the reduction of the sequence
\begin{align*}
\dots{}&\sse{\aslp{i+6}\la}\sse{\zslp{i+4}\la}\sse{-\rslp{i+4}\la}\sse{-\zslp{i+6}\la}
\\
&\sse{\aslp{i+4}\la}\sse{\zslp{i+2}\la}\sse{-\rslp{i+2}\la}\sse{-\zslp{i+4}\la}
\\
&\sse{\aslp{i+2}\la}\sse{\zslp{i}\la}\sse{-\rslp{i}\la}\sse{-\zslp{i+2}\la}
\\
&\sse{\aslp{i}\la}\sse{-\zslp{i}\la}.
\end{align*}
By assumption there are $m<n$ with $m,n\equiv i\ppmod2$ such that $\rslp m\la$ and $\aslp n\la$ are both positive. This means that (letting $s$ denote the sign sequence immediately above) we can choose integers $a_1,\dots,a_t,b_1,\dots,b_t$ with $t=1+\sum_{l\equiv i\ppmod2}\zslp l\la$, such that the hypotheses of \cref{signseq} are satisfied. Hence the reduction of $s$ (which is the reduced $i$-signature of $\mu$) contains fewer than $\epsilon_i\spn\la$ minus signs.
\end{pf}

\subsection{Dimension arguments}\label{degsec}

For some cases in the proof of our main theorem, we employ dimension arguments similar to those used in \cite{mfspin2}. The idea in \cite{mfspin2} is very simple: the $2$-modular reduction $\ol{\spn\la}$ cannot equal $\sid{\la\dblreg}$ if there is another character $\spn\mu$ of smaller degree such that $\ol{\spn\mu}$ also has $\sid{\la\dblreg}$ as a composition factor. Here we extend this idea to show that $\spn\la$ cannot be homogeneous in certain cases, by keeping track of the multiplicity of $\sid{\la\dblreg}$.

Given a $\la\in\cald(n)$, we define $\dim\la$ to be the degree of $\spn\la$, i.e.\ the value of this character at the identity element of $\tsss n$. This is given by the \emph{bar-length formula}, which goes back to Schur \cite{schu}. With our unusual notation $\spn\la$, this formula reads
\[
\deg\spn\la=2^{(n-\len\la)/2}\frac{n!}{\prod_{1\ls i\ls m}\la_i!}\prod_{1\ls i<j\ls m}\frac{\la_i-\la_j}{\la_i+\la_j}.
\]
So $\deg\spn\la$ lies in $\bbn$ if $\ev\la$ is even, or $\sqrt2\bbn$ if $\ev\la$ is odd. Now define the \emph{divided degree}
\[
\ddim\la=\frac{\dim\la}{\dspn_{\la{\la\dblreg}}}=2^{-\ev\la/2}\dim\la.
\]
Now we have the following; this is a simple modification of \cite[Lemma 4.1]{mfspin2}.

\begin{lemma}\label{samerdoub}
Suppose $\la,\mu\in\cald$, with $\la\dblreg=\mu\dblreg$ and $\ddim\la>\ddim\mu$. Then $\spn\la$ is inhomogeneous.
\end{lemma}

Now we give several results which allow us to exploit \cref{samerdoub} in certain situations. We begin with a row-removal lemma.

\begin{lemma}\label{dimrowrem}
Suppose $\la,\mu\in\cald(n)$ with $\mu\doms\la$, and $m$ is an integer with $m>\mu_1$. Define
\begin{align*}
\la^+&=(m,\la_1,\la_2,\dots),\\
\mu^+&=(m,\mu_1,\mu_2,\dots).
\end{align*}
Then
\[
\frac{\ddim{\la^+}}{\ddim{\mu^+}}>\frac{\ddim\la}{\ddim\mu}.
\]
Furthermore, if $\la\dblreg=\mu\dblreg$, then $(\la^+)\dblreg=(\mu^+)\dblreg$.
\end{lemma}

\begin{pf}
A weaker version of the first statement (using degree rather than divided degree) is proved in \cite[Lemma 4.8]{mfspin2}. We consider the ratio of the left-hand side to the right-hand side, which we evaluate using the bar-length formula. Using the fact that $\len{\la^+}-\len{\mu^+}=\len\la-\len\mu$ and $\ev{\la^+}-\ev{\mu^+}=\ev\la-\ev\mu$, we find that this ratio equals
\[
\prod_{i\gs1}\frac{m+\mu_i}{m-\mu_i}\frac{m-\la_i}{m+\la_i}.
\]
This ratio is greater than $1$, by \cite[Lemma 4.7]{mfspin2}.

The second statement is the same as in \cite[Lemma 4.8]{mfspin2}.
\end{pf}

Next we import three results from \cite{mfspin2} comparing divided degrees for specific partitions. We remark that the results we use from \cite{mfspin2} are all stated for degree rather than divided degree, and \cite{mfspin2} does not use our unusual definition of $\spn\la$; but in fact if $\la$ has at most one even part (which is the case for all partitions appearing in the next three lemmas) then $\ddim\la=\deg\lan\la_{(\pm)}\ran$, so the results in \cite{mfspin2} can be imported directly.

\begin{lemmac}{mfspin2}{Lemma 4.2}\label{s42}
Given $a\gs2$, define
\[
\la^a=(4a,4a-3\df5),\quad
\mu^a=(4a+1\df9,4).
\]
Then $(\la^a)\dblreg=(\mu^a)\dblreg$, and $\ddim{\la^a}>\ddim{\mu^a}$.
\end{lemmac}

\begin{lemma}\label{s43}
Given $a\gs1$, define
\[
\la^a=(4a,4a-3\df1),\quad
\mu^a=(4a+1\df5).
\]
Then $(\la^a)\dblreg=(\mu^a)\dblreg$ and $\ddim{\la^a}>\ddim{\mu^a}$.
\end{lemma}

\begin{pf}
This is just the case $m=0$ of \cite[Proposition 4.3]{mfspin2}.
\end{pf}

\begin{lemma}\label{s46}
Given $a\gs1$, define
\[
\la^a=(4a+2,4a-1\df3),\quad
\mu^a=(4a+3\df7,2).
\]
Then $(\la^a)\dblreg=(\mu^a)\dblreg$, and $\ddim{\la^a}>\ddim{\mu^a}$.
\end{lemma}

\begin{pf}
This is the case $m=0$ of \cite[Proposition 4.6]{mfspin2}.
\end{pf}

Now we prove a new result on the same lines.

\begin{lemma}\label{firstdimlem}
Given $a\gs1$, define
\[
\la^a=(4a+2,4a-1\df7,2),\quad\mu^a=(4a+3\df7,1).
\]
Then $\la\dblreg=\mu\dblreg$ and $\ddim\la>\ddim\mu$.
\end{lemma}

\begin{pf}
The first claim is very easy to check. For the second, we use induction on $a$. The base case is a simple check, and the inductive step follows from the claim that
\[
\frac{\ddim{\la^{a+1}}\ddim{\mu^a}}{\ddim{\la^a}\ddim{\mu^{a+1}}}>1.
\]
In fact, we can calculate this ratio directly from the bar-length formula: it is
\[
\frac{(a+1)^2(4a-1)(4a+1)(4a+7)(4a+9)}{a^2(2a+3)(2a+5)(8a+5)(8a+9)}.
\]
To see that this ratio is always greater than $1$, we subtract the denominator from the numerator to get
\[
64a^5+364a^4+544a^3+126a^2-190a-63
\]
which is obviously positive for $a\gs1$.
\end{pf}

\section{Rouquier blocks and \qs partitions}\label{rouqsec}

\subsection{Decomposition numbers for Rouquier blocks}\label{rouqdecsec}

In this section we work with Rouquier blocks of $\tsss n$. These are blocks whose decomposition numbers are relatively well-behaved, and which play an important role in the classification of homogeneous spin characters.

Our main result (\cref{mainrouqdec}) is a formula for the spin decomposition numbers in Rouquier blocks. From this we will be able to deduce which spin characters in Rouquier blocks are homogeneous. We then extend these results to characters labelled by a family of partitions which we call \emph{\qs}.

We begin by recalling the essential definitions and background on Rouquier blocks. Rouquier blocks of symmetric groups have been studied in numerous places, but Rouquier blocks of the double covers are first treated in \cite{mfspin2}. First we fix some notation.

\smallskip
\begin{mdframed}[innerleftmargin=3pt,innerrightmargin=3pt,innertopmargin=3pt,innerbottommargin=3pt,roundcorner=5pt,innermargin=-3pt,outermargin=-3pt]
Throughout \cref{rouqdecsec,rouqhomogsec} we fix a $2$-core $\sigma$, and we let $\tau=\sigma\dbl$ be the corresponding \fbc. For $w\gs0$ we write $\bks w$ for the block of $\tsss{|\sigma|+2w}$ with $2$-core $\sigma$ and weight $w$.
\end{mdframed}

We say that the block $\bks w$ is \emph{Rouquier} if $w\ls\len\sigma+1$. One thing that makes Rouquier blocks easy to understand is the simple description of the partitions labelling the characters.

\begin{propn}\label{rouqptns}
Suppose $0\ls w\ls\len\sigma+1$.
\begin{enumerate}[ref=(\arabic*)]
\item\label{rouqp1}
\cite[Lemma 5.1]{mfspin2} The \trps labelling irreducible Brauer characters in $\bks w$ are precisely the partitions $\sigma+2\mu$, for $\mu\in\calp(w)$.
\item\label{rouqp2}
\cite[Corollaries 5.3 and 5.5]{mfspin2} The \trps labelling spin characters in $\bks w$ are precisely the partitions $\tau+4\alpha\sqcup2\beta$, for $\alpha\in\calp$ and $\beta\in\cald$ with $2|\alpha|+|\beta|=w$.
\end{enumerate}
\end{propn}

So the spin decomposition number problem in Rouquier blocks amounts to finding the decomposition numbers $\decs{(\tau+4\alpha\sqcup2\beta)}{(\sigma+2\mu)}$ for all $\mu,\alpha\in\calp$ and $\beta\in\cald$ with $|\mu|=2|\alpha|+|\beta|\ls\len\sigma+1$.

The corresponding result for decomposition numbers of $\sss n$ in characteristic $2$ is as follows. Recall that for partitions $\la,\mu\in\calp(w)$, $D_{\la\mu}$ denotes the decomposition number for the Schur algebra $\schur w$ (which coincides with the corresponding decomposition number for $\sss n$ if $\mu$ is \treg), and $A_{\la\mu}$ denotes the $(\la,\mu)$-entry of the adjustment matrix.

\begin{thmciting}{\textup{\textbf{\cite[Corollary 2.6]{jmq-1}, \cite[Theorem 132]{turn}.}}\ }\label{lineardecomprouq}
Suppose $w\ls\len\sigma+1$ and $\la,\mu\in\calp(w)$. Then
\[
\odc_{(\sigma+2\la)(\sigma+2\mu)}=\delta_{\la\mu},\qquad
D_{(\sigma+2\la)(\sigma+2\mu)}=D_{\la\mu}
\]
and hence
\[
A_{(\sigma+2\la)(\sigma+2\mu)}=D_{\la\mu}.
\]
\end{thmciting}

Our calculations in Rouquier blocks will be based around projective characters. Recall that a character of $\tsss n$ is projective (in characteristic $2$) if it vanishes on elements of even order. In fact we work with \emph{virtual} projective characters, i.e.\ $\bbz$-linear combinations of projective characters.

Given $\mu\in\cald(n)$, the projective cover of the James module $\jms\mu$ may be lifted to an ordinary representation of $\tsss n$, and we write $\prj\mu$ for the character of this representation; this is called an \emph{indecomposable} projective character, and the characters $\prj\mu$ for $\mu\in\cald(n)$ give a basis for the space of virtual projective characters.

Brauer reciprocity says that $\prj\mu$ is given in terms of irreducible characters by the entries in the column of the decomposition matrix corresponding to $\mu$. With our unusual definition of $\spn\la$ and $\decs\la\mu$, this amounts to the following statement:
\[
\prj\mu=\sum_{\la\in\calp(n)}D_{\la\mu}\ord\la\ \ {+}\sum_{\la\in\cald(n)}\decs\la\mu\spn\la.
\]

Now suppose $w\ls\len\sigma+1$, and consider projective characters in the Rouquier block $\bks w$. By \cref{rouqptns}\ref{rouqp1}, these are linear combinations of the characters $\prj{\sigma+2\mu}$ for $\mu\in\calp(w)$. For each $\mu\in\calp(w)$ we define a virtual projective character $\omega^\mu$ in $\bks w$ by
\begin{align*}
\omega^\mu&=\sum_{\la\in\calp(w)}D\v_{\la\mu}\prj{\sigma+2\la}.
\\
\intertext{Then}
\prj{\sigma+2\mu}&=\sum_{\la\in\calp(w)}D_{\la\mu}\omega^\la
\end{align*}
and \cref{lineardecomprouq} implies that $\omega^\mu$ is the unique virtual projective character in $\bks w$ satisfying
\[
\chm{\omega^\mu}{\ord{\sigma+2\nu}}=\delta_{\mu\nu}
\]
for all $\nu$.

Now we can state our main theorem for Rouquier blocks. (Recall the notation for symmetric functions from \cref{symfnsec}.)

\begin{thm}\label{mainrouqdec}
Suppose $w\ls\len\sigma+1$.  Suppose $\mu,\alpha\in\calp$ and $\beta\in\cald$ with $|\mu|=2|\alpha|+|\beta|=w$. Then
\begin{align*}
\chm{\omega^\mu}{\spn{\tau+4\alpha\sqcup2\beta}}&=2^{\len\beta/2}\sum_{\gamma,\zeta\in\calp}\lan P_\beta,s_\zeta\ran \lrc\mu\gamma\zeta\epstar\gamma\lrc\alpha{\quo\gamma0}{\quo\gamma1}.
\\
\intertext{Hence}
\decs{(\tau+4\alpha\sqcup2\beta)}{(\sigma+2\mu)}&=2^{\len\beta/2}\sum_{\la\in\calp(w)}D_{\la\mu}\sum_{\gamma,\zeta\in\calp}\lan P_\beta,s_\zeta\ran \lrc\la\gamma\zeta\epstar\gamma\lrc\alpha{\quo\gamma0}{\quo\gamma1}.
\end{align*}
\end{thm}

\begin{eg}
Take $w=4$, with $\alpha=(1)$, $\beta=(2)$ and $\mu=(2,1^2)$, and examine the terms in the formula for $\decs{(\tau+4\alpha\sqcup2\beta)}{(\sigma+2\mu)}$. The decomposition matrix for $\schur4$ shows that $D_{\la\mu}=1$ for $\la=(2,1^2)$ and $(1^4)$, and $D_{\la\mu}=0$ otherwise. To get $\lrc{(1)}{\quo\gamma0}{\quo\gamma1}$ non-zero, we need $\quo\gamma0$ and $\quo\gamma1$ to equal $\varnothing$ and $(1)$ in some order. Assuming $\gamma$ also has empty $2$-core (so that $\epstar\gamma\neq0$), we obtain $\gamma=(2)$ or $(1^2)$. The Schur P-function $P_{(2)}$ is easily seen to equal $s_{(2)}+s_{(1^2)}$. So we need to consider pairs of partitions $\gamma,\zeta\in\{(2),(1^2)\}$. For such a pair we get $\lrc{(2,1^2)}\gamma\zeta=1$ as long as at least one of $\gamma$ and $\zeta$ is $(1^2)$, and we get $\lrc{(1^4)}\gamma\zeta=1$ only if they both equal $(1^2)$. So the final summation consists of four terms: three of these (with $\gamma=(1^2)$) contribute a coefficient of $1$, and the other one (with $\gamma=(2)$) gives $-1$. Combining this with the initial coefficient $2^{\len\beta/2}=\sqrt2$, we end up with $\decs{(\tau+4\alpha\sqcup2\beta)}{(\sigma+2\mu)}=2\sqrt2$.
\end{eg}

In fact, we already know a special case of \cref{mainrouqdec} from \cite{mfspin2}.

\begin{propn}\label{mainrouqdecspecial}
\cref{mainrouqdec} holds in the case $\beta=\varnothing$.
\end{propn}

\begin{pf}
Since $P_\varnothing=s_\varnothing$, \cref{mainrouqdec} in the case $\beta=\varnothing$ asserts that
\begin{align*}
\chm{\omega^\mu}{\spn{\tau+4\alpha}}&=\sum_{\gamma}\lrc\mu\gamma\varnothing\epstar\gamma\lrc\alpha{\quo\gamma0}{\quo\gamma1}
\\
&=\epstar\mu\lrc\alpha{\quo\mu0}{\quo\mu1}
.\tag*{($\dagger$)}
\end{align*}
\cite[Theorem 5.14]{mfspin2} says that for any $\mu$,
\[
\decs{(\tau+4\alpha)}{(\sigma+2\mu)}=A_{(\dup\alpha)\mu}.
\]
From above, $\decs{(\tau+4\alpha)}{(\sigma+2\mu)}=\chm{\prj{\sigma+2\mu}}{\spn{\tau+4\alpha}}$, and $\prj{\sigma+2\mu}=\sum_\la D_{\la\mu}\omega^\la$. Hence
\begin{align*}
\chm{\omega^\mu}{\spn{\tau+4\alpha}}&=\sum_\la D\v_{\la\mu}\chm{\prj{\sigma+2\la}}{\spn{\tau+4\alpha}}\\
&=\sum_\la D\v_{\la\mu}A_{(\dup\alpha)\la}\\
&=(AD\v)_{(\dup\alpha)\mu}\\
&=\odc\v_{(\dup\alpha)\mu}.
\end{align*}
By \cref{steincory}, $\odc\v_{(\dup\alpha)\mu}$ equals the right-hand side of ($\dagger$).
\end{pf}

In order to prove \cref{mainrouqdec} by induction, we will need some results on inducing characters in Rouquier blocks. For $a\gs0$ we define $\ffd\bullet a$ to be the function $\ffd0a\ffd1a$ if $\len\sigma$ is odd, or $\ffd1a\ffd0a$ if $\len\sigma$ is even. Given $w\ls\len\sigma+1$, our technique to find characters in the Rouquier block $\bks w$ will be to apply $\ffd\bullet a$ to characters in the block $\bks{w-a}$ (which is also a Rouquier block).

First we show how to apply $\ffd\bullet a$ to a character $\omega^\la$. This result is essentially contained in \cite{mfspin2} (and is effectively a special case of \cite[Lemma 3.1]{ct}).

\begin{propn}\label{ffdomega}
Suppose $\la\in\calp$ with $|\la|+a\ls\len\sigma+1$. Then
\[
\ffd\bullet a\omega^\la=\sum_\mu\lrc\mu\la{(1^a)}\omega^\mu.
\]
\end{propn}

\begin{pf}
The function $\ffd\bullet a$ takes characters in $\bks{\card\la}$ to characters in $\bks{\card\la+a}$. It also takes projective characters to projective characters, so $\ffd\bullet a\omega^\la$ must be a linear combination of the virtual characters $\omega^\mu$ for $\mu\in\calp(\card{\la}+a)$. Since $\omega^\la$ is characterised by $\chm{\omega^\la}{\ord{\sigma+2\mu}}=\delta_{\la\mu}$, we need to show that
\[
\chm{\ffd\bullet a\omega^\la}{\ord{\sigma+2\mu}}=\lrc\mu\la{(1^a)}
\]
for each $\mu$. By definition $\omega^\la$ equals $\ord{\sigma+2\la}$ plus a combination of characters of the form $\ord\epsilon$ with $\epsilon$ $2$-singular, plus a linear combination of spin characters. By \cite[Lemma 5.9]{mfspin2},
\[
\chm{\ffd\bullet a\ord{\sigma+2\la}}{\ord{\sigma+2\mu}}=\lan e_as_\la,s_\mu\ran,
\]
where $e_a$ is the $a$th elementary symmetric function, which coincides with the Schur function $s_{(1^a)}$. So
\[
\chm{\ffd\bullet a\ord{\sigma+2\la}}{\ord{\sigma+2\mu}}=\lrc\mu\la{(1^a)},
\]
while (also from \cite[Lemma 5.9]{mfspin2}) $\chm{\ffd\bullet a\chi}{\ord{\sigma+2\mu}}=0$ for any other irreducible character $\chi$ appearing in $\omega^\la$. The result follows.
\end{pf}

Next we show how to apply $\ffd\bullet a$ to a spin character.

\begin{propn}\label{rouqspinbranch}
Suppose $\alpha,\gamma\in\calp$ and $\beta,\zeta\in\cald$, and $2|\alpha|+|\beta|\ls2|\gamma|+|\zeta|\ls\len\sigma+1$. Let $a=2|\gamma|-2|\alpha|+|\zeta|-|\beta|$. Then
\[
\chm{\ffd\bullet a\spq\tau\alpha\beta}{\spq\tau\gamma\zeta}=
2^{(\len\zeta-\len\beta)/2}\lrcb\gamma\alpha\strp\zeta\beta.
\]
\end{propn}

\begin{pf}
For this proof write $\la=\spq\tau\alpha\beta$ and $\mu=\spq\tau\gamma\zeta$. We assume $\len\sigma$ is even throughout the proof; the other case is the same but with the residues $0$ and $1$ swapped throughout.

In order for the left-hand side of the formula in the \lcnamecref{rouqspinbranch} to be non-zero, we must have $\la\subseteq\mu$ by the spin branching rules. In order for the right-hand side to be non-zero, we must have $\alpha\subseteq\gamma$ and $\beta\subseteq\zeta$, which also gives $\la\subseteq\mu$. So the formula is true (with both sides zero) if $\la\not\subseteq\mu$.

So we assume for the rest of the proof that $\la\subseteq\mu$. Let $s=\len\gamma$. Then we claim that:
\begin{itemize}
\item
up to row $s$, $\la$ agrees with $\tau+4\alpha$ and $\mu$ agrees with $\tau+4\gamma$;
\item
from row $s+1$ on, $\la$ agrees with $\tau\sqcup2\beta$ and $\mu$ agrees with $\tau\sqcup2\zeta$.
\end{itemize}
First we address $\mu$. One of the parts of $\mu$ is $2\zeta_1$, and we claim that this does not occur among the first $s$ parts of $\mu$. The relationship between $\sigma$ and $\tau$ means that $\tau_r=2\len\sigma-4r+3$ for $r\ls\len\tau$. So
\begin{align*}
2\zeta_1&\ls2|\zeta|\\
&\ls2\len\sigma+2-4|\gamma|\\
&\ls2\len\sigma+2-4s\\
&=\tau_s-1.
\end{align*}
Hence the $s$ largest parts of $\mu$ are $\tau_1+4\gamma_1,\dots,\tau_{s}+4\gamma_{s}$. So
\[
(\mu_1,\dots,\mu_s)=(\tau_1+4\gamma_1,\dots,\tau_s+4\gamma_s),\qquad(\mu_{s+1},\mu_{s+2},\dots)=(\tau_{s+1},\tau_{s+2},\dots)\sqcup2\zeta,
\]
which proves our claim for $\mu$. To prove the same claim for $\la$, we first need to show that $2\beta_1<\tau_s$. If this is not true, then $\la$ contains at least $s+1$ parts greater than or equal to $\tau_s$, and in particular $\la_{s+1}\gs\tau_s$. Our assumption that $\la\subseteq\mu$ then means that $\mu_{s+1}\gs\tau_s$; but from above $\mu_{s+1}=\max\{\tau_{s+1},2\zeta_1\}<\tau_s$, a contradiction.

So $2\beta_1<\tau_s$, and hence the first $s$ parts of $\la$ are $\tau_1+4\alpha_1,\dots,\tau_s+4\alpha_s$, while
\[
(\la_{s+1},\la_{s+2},\dots)=(\tau_{s+1}+4\alpha_{s+1},\tau_{s+2}+4\alpha_{s+2},\dots)\sqcup2\beta.
\]
If $\alpha_{s+1}>0$, this gives $\la_{s+1}\gs\tau_s$, but we have just shown that this cannot be the case. So $\len\alpha\ls s$ and
\[
(\la_{s+1},\la_{s+2},\dots)=(\tau_{s+1},\tau_{s+2},\dots)\sqcup2\beta,
\]
so our claim is proved.

Having proved our claim, it is easy to see that to have $\la\subseteq\mu$ we must have $\alpha\subseteq\gamma$ and (using \cref{cupcont}) $\beta\subseteq\zeta$. So we assume this is the case. If $\chm{\ffd1r\ffd0r\spn\la}{\spn\mu}>0$, there must be a \treg partition $\nu$ such that $\la\ads0r\nu\ads1r\mu$. The partition $\nu$ is unique if it exists: it is obtained from $\la$ by adding all the nodes of $\mu\setminus\la$ of \spr $0$. This gives $\nu_r\ls\la_r+2$ and $\mu_r\ls\nu_r+2$ for every $r$, and in particular $\gamma_r\ls\alpha_r+1$ for every $r$; so $\gamma\setminus\alpha$ must be a vertical strip in order for $\nu$ to exist, in which case $\nu_r=\la_r+2(\gamma_r-\alpha_r)$ for $r=1,\dots,s$.

We also claim that $\zeta\setminus\beta$ must be a horizontal strip if $\nu$ exists. The construction of $\nu$ means that for $r=1,\dots,\tau_s$,
\[
\nu'_r=
\begin{cases}
\tau'_r+(2\zeta)'_r&\text{if }r\equiv0,1\ppmod4\\
\tau'_r+(2\beta)'_r&\text{if }r\equiv2,3\ppmod4.
\end{cases}
\]
Since $\nu$ must be a partition, we must also have $\nu'_r\ls\nu'_{r-1}$ for all $r\gs2$. For $r\not\equiv0\ppmod4$ this is immediate, while for $r\equiv0\ppmod4$ this condition says
\[
\tau'_r+(2\zeta)'_r\ls\tau_{r-1}'+(2\beta)'_{r-1}.
\]
The assumption that $\len\sigma$ is even means that $\tau'_r=\tau_{r-1}'-1$ when $r\equiv0\ppmod4$, so our condition becomes
\[
(2\zeta)'_r\ls(2\beta)'_{r-1}+1,
\]
i.e.\
\[
\zeta'_{r/2}\ls\beta'_{r/2}+1
\]
whenever $r\equiv0\ppmod4$. The partition $\nu$ must also be \treg, i.e.\ $\nu'_r\gs\nu'_{r-1}-1$ for all $r\gs2$. This is immediate when $r\not\equiv2\ppmod4$, while for $r\equiv2\ppmod4$ it says
\[
\tau'_r+(2\beta)'_r\gs\tau'_{r-1}+(2\zeta)'_{r-1}-1.
\]
This reduces to
\[
\beta'_{r/2}\gs\zeta'_{r/2}-1
\]
whenever $r\equiv2\ppmod4$. So in order for our \trp $\nu$ to exist we need $\beta'_r\gs\zeta'_r-1$ for all $r$, so that $\zeta\setminus\beta$ is a horizontal strip.

So the formula in the \lcnamecref{rouqspinbranch} is true (with both sides zero) if either $\gamma\setminus\alpha$ is not a vertical strip or $\zeta\setminus\beta$ is not a horizontal strip. So we may assume from now on that $\gamma\setminus\alpha$ is a vertical strip and $\zeta\setminus\beta$ is a horizontal strip. We then need to show that $\chm{\ffd1a\ffd0a\spn\la}{\spn\mu}=2^{(\len\zeta-\len\beta)/2}\strp\zeta\beta$. We do this using two applications of \cref{spinbranchpower}. For this, we just need to count the number of pairs of consecutive columns both containing a node of $\nu\setminus\la$, and the number of pairs of consecutive columns both containing a node of $\mu\setminus\nu$. We start with the nodes of $\nu\setminus\la$. For $r=1,\dots,s$ we have $\nu_r-\la_r\in\{0,2\}$, so the added nodes come in $|\gamma|-|\alpha|$ pairs lying in adjacent columns. In rows $s+1,s+2,\dots$, we add at most one node in each column. We have a pair of added $0$-nodes in adjacent columns for every $r\equiv0\ppmod4$ with $\nu'_r>\la'_r$ and $\nu'_{r+1}>\la'_{r+1}$; writing $r=4t$, this is the same as saying
\[
\tau'_{4t}+(2\zeta)'_{4t}>\tau'_{4t}+(2\beta)'_{4t},\qquad\tau'_{4t+1}+(2\zeta)'_{4t+1}>\tau'_{4t+1}+(2\beta)'_{4t+1},
\]
which is the same as
\[
\zeta'_{2t}>\beta'_{2t},\qquad\zeta'_{2t+1}>\beta'_{2t+1}.
\]
So the number of pairs of nodes added in consecutive columns equals the number of even values of $x$ such that $\zeta\setminus\beta$ contains nodes in columns $x$ and $x+1$. Call this number $X_{\zeta\beta}^{\operatorname{ev}}$. Define $X_{\zeta\beta}^{\operatorname{od}}$ similarly, and observe that $X_{\zeta\beta}^{\operatorname{ev}}+X_{\zeta\beta}^{\operatorname{od}}+N_{\zeta\beta}+\len\zeta-\len\beta=|\zeta|-|\beta|$. Now \cref{spinbranchpower} gives
\begin{align*}
\chm{\ffd0a\spn\la}{\spn\nu}&=2^{(a-\len\nu+\len\la)/2-|\gamma|+|\alpha|-X_{\zeta\beta}^{\operatorname{ev}}}.
\\
\intertext{In a similar way we get}
\chm{\ffd1a\spn\nu}{\spn\mu}&=2^{(a-\len\mu+\len\nu)/2-|\gamma|+|\alpha|-X_{\zeta\beta}^{\operatorname{od}}}.
\\
\intertext{Since there is no other $\xi$ with $\chm{\ffd0a\spn\la}{\spn\xi}\chm{\ffd1a\spn\xi}{\spn\mu}>0$, we obtain}
\chm{\ffd1a\spn\la}{\spn\mu}&=2^{(\len\beta-\len\zeta)/2+|\zeta|-|\beta|-X_{\zeta\beta}^{\operatorname{ev}}-X_{\zeta\beta}^{\operatorname{od}}}\\
&=2^{(\len\zeta-\len\beta)/2+N_{\zeta\beta}}
\end{align*}
as required.
\end{pf}

Now we can use \cref{mainrouqdecspecial} to give the deferred proof of \cref{special2main}.

\begin{pf}[Proof of \cref{special2main}]
If $|\la|>2|\gamma|$ then both sides are zero. Otherwise, let $a=2|\gamma|-|\la|$, and consider the coefficient of $\spn{\tau+4\gamma}$ in $\ffd\bullet a\omega^\la$.

As a consequence of \cref{mainrouqdecspecial,ffdomega}, we obtain
\[
\chm{\ffd\bullet a\omega^\la}{\spn{\tau+4\gamma}}=\sum_\mu\lrc\mu\la{(1^a)}\epstar\mu\lrc\gamma{\quo\mu0}{\quo\mu1},
\]
which is the left-hand side of the formula in \cref{special2main}. On the other hand, using Propositions \ref{rouqptns}\ref{rouqp2}, \ref{mainrouqdecspecial} and \ref{rouqspinbranch}, we obtain
\begin{align*}
\chm{\ffd\bullet a\omega^\la}{\spn{\tau+4\gamma}}&=\sum_\alpha\epstar\la\lrc\alpha{\quo\la0}{\quo\la1}\chm{\ffd\bullet a\spn{\tau+4\alpha}}{\spn{\tau+4\gamma}}
\\
&=\sum_\alpha\epstar\la\lrc\alpha{\quo\la0}{\quo\la1}\lrcb\gamma\alpha
\\
&=\epstar\la\lrcb\gamma{\quo\la0\quo\la1},
\end{align*}
which is the right-hand side. 
\end{pf}

Now we come to the proof of our main result on Rouquier blocks.

\begin{pf}[Proof of \cref{mainrouqdec}]
We just need to prove the first statement; the second statement then follows from the remarks preceding the theorem.

Write $\fs\mu\alpha\beta$ for the right-hand side of the equation in \cref{mainrouqdec}, defining $\fs\mu\alpha\beta=0$ if $|\mu|\neq2|\alpha|+\beta|$. We proceed by induction on $w$, and for fixed $w$ by induction on $\mu$ (with respect to the dominance order). For the case $w=0$, the only irreducible $2$-modular Brauer character in $\bks w$ is $\sid\sigma$, and $\sigma=\tau\dblreg$, so that by \cref{spinreg}
\[
\chm{\omega^\varnothing}{\spn\tau}=\chm{\prj\sigma}{\spn\tau}=1,
\]
in agreement with the theorem.

Now assume $w>0$ and that $\chm{\omega^\nu}{\spq\tau\zeta\eta}=\fs\nu\zeta\eta$ for all $\eta,\zeta$ whenever $|\nu|<w$ or $\nu\domsby\mu$. Suppose the last non-empty column of $\mu$ has length $a$, and let $\mu^-$ denote the partition obtained by removing this last column. By our inductive hypothesis
\[
\chm{\omega^{\mu^-}}{\spq\tau{\alpha^-}{\beta^-}}=\fs{\mu^-}{\alpha^-}{\beta^-}
\]
for all $\alpha^-\in\calp$ and $\beta^-\in\cald$ with $2|\alpha^-|+|\beta^-|=|\mu^-|$. Hence we can compute
\begin{align*}
\chm{\ffd\bullet a\omega^{\mu^-}}{\spq\tau\alpha\beta}&=\sum_{\substack{\alpha^-\in\calp\\\beta^-\in\cald}}\fs{\mu^-}{\alpha^-}{\beta^-}\chm{\ffd\bullet a\spq\tau{\alpha^-}{\beta^-}}{\spq\tau\alpha\beta}\\
&=\sum_{\substack{\alpha^-\in\calp\\\beta^-\in\cald}}2^{\len{\beta^-}/2}\sum_{\xi,\eta\in\calp}\lan P_{\beta^-},s_\eta\ran\lrc{\mu^-}\xi\eta\epstar\xi\lrc{\alpha^-}{\quo\xi0}{\quo\xi1}
2^{(\len\beta-\len{\beta^-})/2}\strp\beta{\beta^-}\lrcb\alpha{\alpha^-}\\
&=2^{\len\beta)/2}\sum_{\substack{\beta^-\in\cald\\\xi,\eta\in\calp}}\lan P_{\beta^-},s_\eta\ran\lrc{\mu^-}\xi\eta\epstar\xi\lrcb\alpha{\quo\xi0\quo\xi1}
\strp\beta{\beta^-}\\
&=2^{\len\beta)/2}\sum_{\nu,\xi,\eta\in\calp}\lrcb\nu{\mu^-}\ip{P_\beta}{s_\eta}\epstar\xi\lrc\alpha{\quo\xi0}{\quo\xi1}\lrc\nu\xi\eta\tag*{by \cref{mainsym}}\\
&=\sum_{\nu\in\calp}\lrcb\nu{\mu^-}\fs\nu\alpha\beta.
\end{align*}
Since $\fs\nu\alpha\beta=0$ unless $|\nu|=|\mu|$, we can write this as
\[
\sum_{\nu\in\calp(w)}\lrc\nu{\mu^-}{(1^a)}\fs\nu\alpha\beta.
\]

By \cref{ffdomega} we also have
\[
\chm{\ffd\bullet a\omega^{\mu^-}}{\spq\tau\alpha\beta}=\sum_{\nu\in\calp(w)}\lrc\nu{\mu^-}{(1^a)}\chm{\omega^\nu}{\spq\tau\alpha\beta},
\]
so that
\[
\sum_{\nu\in\calp(w)}\lrc\nu{\mu^-}{(1^a)}\chm{\omega^\nu}{\spq\tau\alpha\beta}=\sum_{\nu\in\calp(w)}\lrc\nu{\mu^-}{(1^a)}\fs\nu\alpha\beta.
\]
Since $\mu^-$ is obtained by removing the last column from $\mu$, every partition $\nu$ with $\lrc\nu{\mu^-}{(1^a)}>0$ satisfies $\nu\domby\mu$. The inductive hypothesis gives $\chm{\omega^\nu}{\spq\tau\alpha\beta}=\fs\nu\alpha\beta$ for $\nu\domsby\mu$, and hence $\chm{\omega^\mu}{\spq\tau\alpha\beta}=\fs\mu\alpha\beta$ as well.
\end{pf}

\subsection{Homogeneous spin characters in Rouquier blocks}\label{rouqhomogsec}

In this section we use \cref{mainrouqdec} to classify homogeneous spin characters in Rouquier blocks. Our main result here is as follows.

\begin{thm}\label{rouqhomog}
Suppose $\alpha\in\calp$ and $\beta\in\cald$ with $2|\alpha|+|\beta|\ls\len\sigma+1$. Then $\spq\tau\alpha\beta$ is homogeneous \iff $\beta$ is a $2$-core and $\alpha$ is $2$-Carter.
\end{thm}

\begin{mdframed}[innerleftmargin=3pt,innerrightmargin=3pt,innertopmargin=3pt,innerbottommargin=3pt,roundcorner=5pt,innermargin=-3pt,outermargin=-3pt]
For the rest of \cref{rouqhomogsec} we fix $\alpha\in\calp$ and $\beta\in\cald$ with $2|\alpha|+|\beta|\ls\len\sigma+1$, and we write $\la=\spq\tau\alpha\beta$.
\end{mdframed}

We begin with the case where $\beta$ is a $2$-core.

\begin{propn}\label{2corerow}
Suppose $\beta$ is a $2$-core. Then for $\mu\in\calp(w)$,
\[
\decs{\la}{(\sigma+2\mu)}=2^{\len\beta/2}A_{(\dup\alpha\sqcup\beta)\mu}.
\]
Hence $\spn{\la}$ is homogeneous \iff $\alpha$ is $2$-Carter.
\end{propn}

\begin{pf}
Since $\beta$ is a $2$-core we have $P_\beta=s_\beta$ by \cref{p2core}. We also have $\odc\v_{\beta\kappa}=\delta_{\beta\kappa}$ for any $\kappa$, because $\weyl\beta$ lies in a block of $\qschur{\card\beta}$ by itself. So (with variables re-labelled) the first formula in \cref{steinberggeneral} becomes
\[
\odc\v_{(\dup\alpha\sqcup\beta)\nu}=\sum_{\gamma\in\calp}\epstar\gamma\lrc\alpha{\quo\gamma0}{\quo\gamma1}\lrc\nu\beta\gamma.
\]
Now the formula in \cref{mainrouqdec} becomes
\begin{align*}
\decs{(\tau+4\alpha\sqcup2\beta)}{(\sigma+2\mu)}&=2^{\len\beta/2}\sum_{\nu\in\calp(w)}D_{\nu\mu}\sum_{\gamma\in\calp}\lrc\nu\gamma\beta\epstar\gamma\lrc\alpha{\quo\gamma0}{\quo\gamma1}
\\
&=2^{\len\beta/2}\sum_{\nu\in\calp(w)}\odc_{(\dup\alpha\sqcup\beta)\nu}\v D_{\nu\mu}
\\
&=2^{\len\beta/2}A_{(\dup\alpha\sqcup\beta)\mu}.
\end{align*}
By \cref{steincory}, $A_{(\dup\alpha\sqcup\beta)\mu}$ equals $D_{\alpha\gamma}$ if $\mu$ has the form $\dup\gamma\sqcup\beta$, and $0$ otherwise. So in order for $\spn{\la}$ to be homogeneous, there must be a unique $\gamma$ with $D_{\alpha\gamma}>0$; in other words, $\weyl\alpha$ is homogeneous. By \cref{dndom,irredweyl}, this is the same as saying that $\alpha$ is $2$-Carter.
\end{pf}

Now we come to the case where $\beta$ is not a $2$-core. 

\begin{lemma}\label{zetacover}
Suppose $\eta\in\calp$ is obtained from $\beta$ by moving a node from row $r$ to row $r+1$ for some $r$, and let
\[
\mu=\dup\alpha\sqcup\beta',\qquad\nu=\dup\alpha\sqcup\eta'.
\]
Then $\chm{\omega^\mu}{\spn\la}>0$ and $\chm{\omega^\nu}{\spn\la}>0$.
\end{lemma}

\begin{pf}
By \cref{mainrouqdec}
\[
\chm{\omega^\xi}{\spn{\tau+4\alpha\sqcup2\beta}}=2^{\len\beta/2}\sum_{\gamma,\zeta\in\calp}\lan P_\beta,s_\zeta\ran\lrc\xi\gamma\zeta\epstar\gamma\lrc\alpha{\quo\gamma0}{\quo\gamma1}
\]
for any $\xi$. Suppose we have $\gamma,\zeta$ such that the summand $\lan P_\beta,s_\zeta\ran\lrc\xi\gamma\zeta\epstar\gamma\lrc\alpha{\quo\gamma0}{\quo\gamma1}$ is non-zero. In order to have $\lan P_\beta,s_\zeta\ran\neq0$, we need $\lan P_\beta,s_{\zeta'}\ran\neq0$ by \cref{pomeg}, and hence $\beta\dom\zeta'$ by \cref{schurpdom}. In order to have $\lrc\xi\gamma\zeta\lrc\alpha{\quo\gamma0}{\quo\gamma1}\neq0$, we need $\xi\dom\gamma\sqcup\zeta$ and $\alpha\domby\quo\gamma0+\quo\gamma1$ by \cref{basiclr}. Using \cref{domsqcup,empcore} as well, we find
\[
\xi\dom\gamma\sqcup\zeta\dom(\quo\gamma0+\quo\gamma1)\sqcup(\quo\gamma0+\quo\gamma1)\sqcup\zeta\dom\alpha\sqcup\alpha\sqcup\zeta\dom\alpha\sqcup\alpha\sqcup\beta'.\tag*{($\dagger$)}
\]
In the case $\xi=\mu$, we must therefore have equality throughout, so that $\gamma=\dup\alpha$ and $\zeta=\beta'$. We then obtain
\begin{itemize}
\item
$\lan P_\beta,s_\zeta\ran=1$ by \cref{pomeg,schurpdom};
\item
$\lrc\mu\gamma\zeta\lrc\alpha{\quo\gamma0}{\quo\gamma1}=1$ by \cref{basiclr}, because $\mu=\gamma\sqcup\zeta$ and $\alpha=\quo\gamma0+\quo\gamma1$;
\item
$\epstar\gamma=1$ because all the columns of $\gamma$ have even length.
\end{itemize}
So we get $\chm{\omega^\mu}{\spn\la}>0$.

Now consider $(\dagger)$ with $\xi=\nu$. The relationship between $\beta$ and $\eta$ means that $\alpha\sqcup\alpha\sqcup\beta'$ is obtained from $\alpha\sqcup\alpha\sqcup\eta'$ by moving a node from column $r+1$ to column $r$, so $\alpha\sqcup\alpha\sqcup\eta'$ covers $\alpha\sqcup\alpha\sqcup\beta'$ in the dominance order, by \cref{domcover}. Hence we must have equality in three of the four inequalities in ($\dagger$). This gives at most four possible cases.
\begin{enumerate}
\item
$\alpha\sqcup\alpha\sqcup\eta'=\gamma\sqcup\zeta=(\quo\gamma0+\quo\gamma1)\sqcup(\quo\gamma0+\quo\gamma1)\sqcup\zeta=\alpha\sqcup\alpha\sqcup\zeta$.

These equalities are satisfied \iff $\gamma=\alpha\sqcup\alpha$ and $\zeta=\eta'$. And now we have:
\begin{itemize}
\item
$\lan P_\beta,s_\zeta\ran=1$ by \cref{pomeg,schurpcover};
\item
$\lrc\nu\gamma\zeta\lrc\alpha{\quo\gamma0}{\quo\gamma1}=1$ by \cref{basiclr}, because $\nu=\gamma\sqcup\zeta$ and $\alpha=\quo\gamma0+\quo\gamma1$;
\item
$\epstar\gamma=1$ because all the columns of $\gamma$ have even length.
\end{itemize}
Hence $\lan P_\beta,s_\zeta\ran \lrc\nu\gamma\zeta\epstar\gamma\lrc\alpha{\quo\gamma0}{\quo\gamma1}=1$ in this case.

\item
$\alpha\sqcup\alpha\sqcup\eta'=\gamma\sqcup\zeta=(\quo\gamma0+\quo\gamma1)\sqcup(\quo\gamma0+\quo\gamma1)\sqcup\zeta$ and $\alpha\sqcup\alpha\sqcup\zeta=\alpha\sqcup\alpha\sqcup\beta'$.

Here the second equality tells us that all the columns of $\gamma$ have even length, and then the first equality implies that $\eta_s\equiv\zeta'_s\ppmod2$ for all $s$. But the final equality gives $\zeta'=\beta$, a contradiction.

\item
$\alpha\sqcup\alpha\sqcup\eta'=\gamma\sqcup\zeta$ and $(\quo\gamma0+\quo\gamma1)\sqcup(\quo\gamma0+\quo\gamma1)\sqcup\zeta=\alpha\sqcup\alpha\sqcup\zeta=\alpha\sqcup\alpha\sqcup\beta'$.

Here the last equality gives $\zeta=\beta'$, and the first equality then gives $\dup\alpha\sqcup\eta'=\gamma\sqcup\beta'$, so $\gamma$ is obtained from $\alpha\sqcup\alpha$ by moving a node from column $r$ to column $r+1$. If $\alpha'_r=\alpha'_{r+1}$ then there is no such $\gamma$, so assume $\alpha'_r>\alpha'_{r+1}$. Now \cref{duplem} says that $\alpha=\quo\gamma0+\quo\gamma1$ (so that the second equality is satisfied) and:
\begin{itemize}
\item
$\lan P_\beta,s_\zeta\ran=1$ by \cref{pomeg,schurpdom};
\item
$\lrc\nu\gamma\zeta\lrc\alpha{\quo\gamma0}{\quo\gamma1}=1$ by \cref{basiclr}, because $\nu=\gamma\sqcup\zeta$ and $\alpha=\quo\gamma0+\quo\gamma1$;
\item
$\epstar\gamma=-1$ by \cref{duplem}.
\end{itemize}

Hence $\lan P_\beta,s_\zeta\ran \lrc\nu\gamma\zeta\epstar\gamma\lrc\alpha{\quo\gamma0}{\quo\gamma1}=-1$ in this case.
\item
$\gamma\sqcup\zeta=(\quo\gamma0+\quo\gamma1)\sqcup(\quo\gamma0+\quo\gamma1)\sqcup\zeta=\alpha\sqcup\alpha\sqcup\zeta=\alpha\sqcup\alpha\sqcup\beta'$.

In this case (using \cref{empcore}) the equalities are satisfied \iff $\gamma=\alpha\sqcup\alpha$ and $\zeta=\beta'$. So:
\begin{itemize}
\item
$\lan P_\beta,s_\zeta\ran=1$ by \cref{pomeg,schurpdom};
\item
$\epstar\gamma=1$ because the columns of $\gamma$ have even length;
\item
$\lrc\alpha{\quo\gamma0}{\quo\gamma1}=1$ by \cref{basiclr}, because $\alpha=\quo\gamma0+\quo\gamma1$.
\end{itemize}
Under the additional assumption that $\alpha_r'>\alpha_{r+1}'$, we claim that $\lrc\nu\gamma\zeta>0$; that is, $\lrc{\alpha\sqcup\alpha\sqcup\eta'}{(\alpha\sqcup\alpha)}{\beta'}>0$. Replacing partitions with their conjugates, we just need to show that $\lrc{2\alpha'+\eta}{(2\alpha')}{\beta}>0$, and this follows from \cref{lrlem}.

So in this case we have $\lan P_\beta,s_\zeta\ran \lrc\nu\gamma\zeta\epstar\gamma\lrc\alpha{\quo\gamma0}{\quo\gamma1}\gs0$, and the inequality is strict if $\alpha'_r>\alpha'_{r+1}$.
\end{enumerate}
Summing over the four cases, we conclude that $\chm{\omega^\nu}{\spn{\tau+4\alpha\sqcup2\beta}}>0$, as required.
\end{pf}

\begin{pf}[Proof of \cref{rouqhomog}]
The case where $\beta$ is a $2$-core is addressed in \cref{2corerow}. So assume $\beta$ is not a $2$-core. Then there is $r\gs1$ such that $\beta_r\gs\beta_{r+1}+2$. Let $\eta$ be the partition obtained from $\beta$ by moving a node from row $r$ to row $r+1$, and define $\mu,\nu$ as in \cref{zetacover}. Then we claim that $\decs\la{(\sigma+2\mu)}$ and $\decs\la{(\sigma+2\nu)}$ are both positive, so that $\spn\la$ is inhomogeneous.

For any $\pi$,
\begin{align*}
\decs\la{(\sigma+2\pi)}&=\chm{\prj{\sigma+2\pi}}{\spn\la}\\
&=\sum_\xi D_{\xi\pi}\chm{\omega^\xi}{\spn\la}.
\end{align*}
Now $D_{\xi\pi}>0$ only if $\pi\dom\xi$. On the other hand, the argument in the first paragraph of the proof of \cref{zetacover} shows that $\chm{\omega^\xi}{\spn\la}$ is non-zero only if $\xi\dom\mu$. So (using \cref{dndom} as well) we can write
\[
\decs\la{(\sigma+2\pi)}=\sum_{\pi\dom\xi\dom\mu}D_{\xi\pi}\chm{\omega^\xi}{\spn\la}.\tag*{($\ddagger$)}
\]
In the particular case $\pi=\mu$, this gives
\[
\decs\la{(\sigma+2\mu)}=\chm{\omega^\mu}{\spn\la},
\]
and this is positive by \cref{zetacover}.

In the case $\pi=\nu$, the fact that $\nu$ covers $\mu$ in the dominance order gives
\[
\decs\la{(\sigma+2\nu)}=\chm{\omega^\nu}{\spn\la}+D_{\mu\nu}\chm{\omega^\mu}{\spn\la}.
\]
By \cref{zetacover} the first term is strictly positive, and the second is non-negative, so ${\decs\la{(\sigma+2\nu)}>0}$.
\end{pf}

\subsection{Separated partitions}

In this subsection we extend the results of \cref{rouqdecsec} beyond Rouquier blocks, to a family of \trps we call \emph{\qs}. These were addressed in \cite{mfspin2}, but here we introduce a more general definition of \qs partition, where we allow arbitrarily many even parts.

Suppose $\la\in\cald$ and $i\in\{1,3\}$. Say that $\la$ is \emph{$i$-\qs} if
\begin{itemize}
\item
all the odd parts of $\la$ are congruent to $i$ modulo $4$, and
\item
if any $\la_r$ is even, then $\la$ includes all positive integers less than $\la_r$ which are congruent to $i$ modulo $4$.
\end{itemize}

Say that $\la$ is \emph{\qs} if it either $1$- or $3$-\qs. \Qs partitions can alternatively be characterised as follows.

\begin{lemma}\label{sepalt}
Suppose $\la\in\cald$. Then $\la$ is \qs \iff it can be written in the form $\tau+4\alpha\sqcup2\beta$, where $\tau$ is a \fbc, $\alpha\in\calp$ and $\beta\in\cald$ with $\tau_1\gs4\len\alpha+2\beta_1-3$.
\end{lemma}

\begin{pf}
Suppose $\la$ is $i$-\qs. Write $\la=\kappa\sqcup2\beta$, where all the non-zero parts of $\kappa$ are odd. By assumption the non-zero parts of $\kappa$ are all congruent to $i$ modulo $4$, so we can write $\kappa=\tau+4\alpha$, where $\tau$ is a \fbc and $\len\alpha\ls\len\tau$.

If $\tau=\varnothing$, then $\la$ has no odd parts, so the \qs condition means that $2\beta_1\ls2$, as required. On the other hand, if $\tau\neq\varnothing$, then the first positive integer congruent to $i$ modulo $4$ which is not contained in $\la$ is $\tau_1+4-4\len\alpha$, so the \qs condition implies that $2\beta_1<\tau_1+4-4\len\alpha$, which is the desired inequality.

Conversely, suppose $\la$ can be written as in the \lcnamecref{sepalt}. Let $i$ be the common residue modulo $4$ of the odd parts of $\tau$, setting $i=3$ if $\tau=\varnothing$. Then all the odd parts of $\la$ are congruent to $i$ modulo $4$. To check the second condition in the definition of an $i$-\qs partition, we need to show that $\la$ contains all positive integers congruent to $i$ modulo $4$ which are less than $2\beta_1$. The first positive integer congruent to $i$ modulo $4$ which is not contained in $\tau+4\alpha$ is $\tau_1+4-4\len\alpha$ if $\tau\neq\varnothing$, or $3$ if $\tau=\varnothing$. The inequality $\tau_1\gs4\len\alpha+2\beta_1-3$ implies that $2\beta_1$ is less than this, so that $\la$ does contain all positive integers congruent to $i$ modulo $4$ and less than~$2\beta_1$.
\end{pf}

It is clear that for a \qs \trp the expression $\tau+4\alpha\sqcup2\beta$ in \cref{sepalt} is unique. Now we can state our main result for \qs partitions, which extends \cref{rouqhomog}.

\begin{thm}\label{mainqs}
Suppose $\la$ is a \qs \trp, and write $\la=\tau+4\alpha\sqcup2\beta$ with $\tau$ a \fbc. Then $\spn\la$ is homogeneous \iff $\beta$ is a $2$-core and $\alpha$ is $2$-Carter.
\end{thm}

The proof of \cref{mainqs} is a downwards induction, using \cref{rouqhomog} as an initial case. The set-up for the inductive step generalises the results of \cite[Section 5.7]{mfspin2} in a straightforward way.

\begin{lemma}\label{qsind}
Suppose $\la$ is a \qs \trp, and write $\la=\tau+4\alpha\sqcup2\beta$, where $\tau$ is a \fbc. Let $i$ be the common \spr of the \spams of $\tau$ and let $\upsilon$ be the \fbc obtained by adding all the $i$-\spams to $\tau$. Let $\mu=\upsilon+4\alpha\sqcup2\beta$. Then $\mu$ is \qs, and
\[
\fmx i\spn\la=\spn\mu,\qquad\emx i\spn\mu=\spn\la,
\]
so $\spn\la$ is homogeneous \iff $\spn\mu$ is.
\end{lemma}

\begin{pf}
The fact that $\mu$ is \qs is immediate from \cref{sepalt}. For the statements about induction and restriction, we assume for ease of notation that $i=0$, though the other case is very similar. Let $l=\len\tau$. Observe that $\la$ has $0$-\spams in column $1$ and in columns $\tau_r+4\alpha_r+1$ and $\tau_r+4\alpha_r+2$ for each $r=1,\dots,l$. Adding all these nodes yields $\mu$, so by \cref{spinbranchpower}
\[
\fmx i\spn\la=\ffd i{2l+1}\spn\la=\spn\mu.
\]
In a similar way we get $\emx i\spn\mu=\spn\la$.

For the final statement, suppose $\spn\mu$ is homogeneous. Then by \cref{spinreg} $\mspn\mu$ is a scalar multiple of $\sid{\mu\dblreg}$. Hence $\epsilon_i\sid{\mu\dblreg}=\epsilon_i\spn\mu$, so that $\mspn\la=\emx i\mspn\mu$ is a scalar multiple of $\emx i\sid{\mu\dblreg}$. By \cref{modbranch}, $\emx i\sid{\mu\dblreg}$ is an irreducible Brauer character, so $\spn\la$ is homogeneous. The opposite implication is proved in the same way.
\end{pf}

Now \cref{mainqs} follows by downwards induction starting from \cref{rouqhomog} in exactly the same way as \cite[Proposition 5.27]{mfspin2} is proved.

\section{Proof of the main theorem}\label{mainproofsec}

In this section we complete the proof of \cref{mainthm}. We begin with the ``if'' part, which is essentially done.

\begin{pf}[Proof of \cref{mainthm} (``if'' part)]
If $\la$ has the form $\tau+4\alpha\sqcup(4,2)$ with $\len\tau>\len\alpha$ and $\alpha$ a $2$-Carter partition, then by \cref{sepalt} $\la$ is \qs, so the result follows from \cref{mainqs}. On the other hand, if $\la$ is one of the partitions $(4,3,2,1)$, $(5,4,3,2)$, $(5,4,3,2,1)$, $(7,4,3,2,1)$, then we can check directly using the Modular Atlas Homepage \cite{modatl} that $\spn\la$ is homogeneous. For $n\ls17$, this page gives decomposition numbers for the symmetric group $\sss n$, and also the degrees of the irreducible Brauer characters. One can then check using the bar-length formula that $\ddim\la=\deg\sid{\la\dblreg}$, so that by \cref{spinreg} $\mspn\la=2^{\ev\la/2}\sid{\la\dblreg}$. For these partitions, we get the following values.
\pfendeqn{
\begin{array}{ccc}
\hline
\la&\la\dblreg&\ddim\la=\deg\sid{\la\dblreg}
\\
\hline
(4,3,2,1)&(7,3)&48
\\
(5,4,3,2)&(8,5,1)&4576
\\
(5,4,3,2,1)&(9,5,1)&4576
\\
(7,4,3,2,1)&(9,5,2,1)&339456
\\
\hline
\end{array}}
\end{pf}

It remains to check the ``only if'' part of \cref{mainthm}, by showing that if $\la$ is a \trp that has two non-zero even parts but is not one of the partitions in \cref{mainthm} then $\spn\la$ is inhomogeneous. We start by singling out seven partitions which cannot be dealt with by any of our main inductive arguments. Let
\[
\scrr=\{(8,4),(8,3,2,1),(12,3,2,1),(13,4,3,2,1),(11,5,4,3,2),(15,5,4,3,2),(19,11,5,4,3,2)\}.
\]

\begin{propn}\label{scrrred}
Suppose $\la\in\scrr$. Then $\spn\la$ is inhomogeneous.
\end{propn}

\begin{pf}
In each case we can find a composition factor of $\spn\la$ other than $\sid{\la\dblreg}$ using a Fock space calculation.

As explained in \cref{decnosec}, we can equate $\mspn\la$ with the $2$-modular reduction of a linear combination of ordinary irreducible characters. That is, we can write $\mspn\la=\sum_\xi a_\xi\modr{\ord\xi}$ for some coefficients $a_\xi\in\bbc$. Hence for any $\mu$ we have $\dspn_{\la\mu}=\sum_\xi a_\xi D_{\xi\mu}$. The values $D_{\xi\mu}$ are not known in all the cases required for the present proposition, but we can exploit the theory of adjustment matrices. Recall from \cref{adjmat} that $D_{\xi\mu}=\sum_\nu\odc_{\xi\nu}A_{\nu\mu}$, where $A$ is the adjustment matrix for the Schur algebra. Hence if we let $C_{\la\mu}=\sum_\nu\dspn_{\la\nu}A\v_{\nu\mu}$, then $C_{\la\mu}=\sum_\xi a_\xi\odc_{\xi\mu}$. The values $\odc_{\xi\mu}$ can be computed using the LLT algorithm \cite{ari,llt}, and hence we can compute $C_{\la\mu}$ for any $\mu$.  (We do not give details of this algorithm here, since we are only applying it in a few cases.)

Now suppose we can find $\mu\in\cald(n)$ with $\mu\neq\la\dblreg$ such that $C_{\la\mu}>0$ while $C_{\la\nu}\gs0$ for all $\nu\domby\mu$. Since $A$ has non-negative entries, with $A_{\mu\mu}=1$ and $A_{\nu\mu}>0$ only for $\nu\domby\mu$, we obtain
\[
\decs\la\mu=\sum_{\nu\in\cald(n)}C_{\la\nu}A_{\nu\mu}>0,
\]
so that $\spn\la$ is inhomogeneous.

For each $\la\in\scrr$ we carry out this calculation using the LLT algorithm in GAP \cite{gap}, and it turns out that in each case we can find a partition $\mu$ as above. We give the value of $\mu$ in each case.
\pfendeqn{
\begin{array}{ccc}
\hline
\la&\la\dblreg&\mu\\
\hline
(8,4)&(5,4,2,1)&(7,5)\\
(8,3,2,1)&(7,4,2,1)&(9,5)\\
(12,3,2,1)&(7,6,4,1)&(9,5,3,1)\\
(13,4,3,2,1)&(9,6,5,3)&(9,7,5,2)\\
(11,5,4,3,2)&(10,7,4,3,1)&(10,7,5,3)\\
(15,5,4,3,2)&(10,7,6,5,1)&(10,9,7,3)\\
(19,11,5,4,3,2)&(12,9,8,7,4,3,1)&(12,9,8,7,5,3)\\
\hline
\end{array}}
\end{pf}

Now we come to the inductive part of the proof of \cref{mainthm}. This relies heavily on the next \lcnamecref{atmostev}.

\begin{lemma}\label{atmostev}
Suppose $\la\in\cald$ and $\spn\la$ is homogeneous, and take $i\in\{0,1\}$. Then:
\begin{enumerate}[ref=(\arabic*)]
\item\label{ladnho}
$\spn{\la\dn i}$ is homogeneous;
\item\label{ladnev}
$\ev{\la\dn i}\ls\ev\la$;
\item\label{epsid}
$\epsilon_i\sid{\la\dblreg}=\epsilon_i\spn\la$;
\item\label{eesid}
$\emx i(\sid{\la\dblreg})=\sid{(\la\dn i)\dblreg}$.
\end{enumerate}
\end{lemma}

\begin{pf}
Let $\mu=\la\dn i$, and let $a$ denote the number of $i$-\sprms of $\la$. Since $\spn\la$ is homogeneous, $\mspn\la$ equals $2^{\ev\la/2}\sid{\la\dblreg}$ by \cref{spinreg}. Hence $\epsilon_i\spn\la=\epsilon_i\sid{\la\dblreg}=a$, and $\emx i{\mspn\la}=2^{\ev\la/2}\emx i{\sid{\la\dblreg}}$, with $\emx i{\sid{\la\dblreg}}$ being an irreducible Brauer character by \cref{modbranch}.

On the other hand, by \cref{spinbranchpower},
\[
\emx i{\spn\la}=\eed ia\spn\la=2^{(a-\len\la+\len\mu)/2-c}\spn\mu,\tag*{(\textasteriskcentered)}
\]
where $c$ is the number of pairs of nodes of $\la\setminus\mu$ lying in consecutive columns. So
\[
\mspn\mu=2^{(\ev\la-a+\len\la-\len\mu)/2+c}\emx i{\sid{\la\dblreg}}
\]
is homogeneous. Now applying \cref{spinreg} to $\mu$ we deduce that $\emx i{\sid{\la\dblreg}}=\sid{\mu\dblreg}$, and that
\[
\frac{\ev\la-a+\len\la-\len\mu}2+c=\frac{\ev\mu}2,
\]
so
\[
\ev\la-\ev\mu=a-\len\la+\len\mu-2c,
\]
which is clearly non-negative, giving \ref{ladnev}.
\end{pf}

An alternative version of \cref{atmostev} holds for induction, replacing $\dn i$, $\epsilon_i$, $\ee i$ with $\up i$, $\phi_i$, $\ff i$. We will refer to this version also as \cref{atmostev}.

With this \lcnamecref{atmostev} in mind, we introduce some notation. For $i=0,1$ we let $\ird_i$ denote the set of all \trps $\la$ for which $\ev\la=i$ and $\spn\la$ is homogeneous; these are given by \cite[Theorem 3.3]{mfspin2}, and are simply the partitions in \cref{mainthm01} together with the partitions $(4b)$ for $b\gs1$. We also let $\ird_2$ denote the set of \trps with exactly two even parts which we claim label homogeneous spin characters; these are precisely the partitions appearing in \cref{mainthm}. Let $\ird=\ird_0\cup\ird_1\cup\ird_2$.

Given a \trp $\la$ and $i\in\{0,1\}$, recall that we write $\noregdown i\la$ to mean that $\epsilon_i\sid{\la\dblreg}<\epsilon_i\spn\la$; by \cref{atmostev}, if $\noregdown i\la$ then $\spn\la$ is inhomogeneous.

In order to prove the ``only if'' part of \cref{mainthm} by induction, we take a \trp $\la$, and we can assume (in view of \cref{atmostev}) that for $i=0,1$ either $\la\dn i=\la$ or $\la\dn i\in\ird$. We set out our assumptions for easy reference.

\smallskip
\needspace{10em}
\begin{mdframed}[innerleftmargin=3pt,innerrightmargin=3pt,innertopmargin=3pt,innerbottommargin=3pt,roundcorner=5pt,innermargin=-3pt,outermargin=-3pt]
\noindent\textbf{Assumptions in force for the rest of \cref{mainproofsec}:}

$\la$ is a \trp with $\ev\la=2$. $\la$ is not \qs and does not lie in $\ird_2$ or $\scrr$. For $i=0,1$, either $\rs\la i\in\ird$ or $\rs\la i=\la$.
\end{mdframed}

We need one more item of notation: for any $r\gs1$, we write $\delta_r$ for the composition which has a $1$ in position $r$ and $0$s everywhere else.

First we extract more specific information about $\la\dn0$ and $\la\dn1$ by examining the \sprs of the \spre and \spams of partitions in $\ird$. Suppose that $i\in\{0,1\}$ with $\la\dn i\in\ird$. We consider the three possibilities for $\ev{\la\dn i}$.
\begin{description}
\item[$\la\dn i\in\ird_0$]\indent\\
In this case $\la\dn i$ has the form $\tau+4\alpha$, where $\tau$ is a \fbc with $\tau_1\equiv2i-1\ppmod 4$, and $\alpha$ is a $2$-Carter partition with $\len\alpha\ls\len\tau$.
\item[$\la\dn i\in\ird_1$]\indent\\
In this case $\la\dn i$ cannot be of the form $(2b)$, since there is no way to add nodes of \spr $i$ to $(2b)$ to obtain a partition with two even positive parts. $\la\dn i$ cannot have the form $(4b-2,1)$, since this partition has \sprms with both residues. Finally, one can check that $\la\dn i$ cannot equal $(3,2,1)$.

So $\la\dn i$ has the form $\tau+4\alpha\sqcup(2)$, where $\tau$ is a \fbc with $\tau_1\equiv2i-1\ppmod4$, and $\alpha$ is a $2$-Carter partition with $\len\alpha\ls\len\tau-i$.
\item[$\la\dn i\in\ird_2$]\indent\\
In this case one can check that (given our other assumptions) $\la\dn i$ cannot be any of $(4,3,2,1)$, $(5,4,3,2)$, $(5,4,3,2,1)$, $(7,4,3,2,1)$. So $\la\dn i$ has the form $\tau+4\alpha\sqcup(4,2)$, where $\tau$ is a \fbc with $\tau_1\equiv2i-1\ppmod4$ and $\alpha$ is a $2$-Carter partition with $\len\alpha<\len\tau$.
\end{description}

Now we consider the various possible cases for $\la\dn0$ and $\la\dn1$.

\subsection*{\nextcase$\la\dn0=\la$, $\la\dn1\in\ird_0$}

In this case we write
\[
\la\dn1=(4m-3\df1)+4\alpha
\]
where $\alpha$ is a $2$-Carter partition with $\len\alpha\ls m$. We reconstruct $\la$ from this partition by adding nodes of \spr $1$, and since $\la\dn0=\la$ we need to add enough  of these nodes to ensure that $\la$ has no $0$-\sprms. Together with the fact that $\ev\la=2$, this means that we add two nodes in each row from $1$ to $m$ except for two rows where we only add one. So
\[
\la=(4m-1\df3)+4\alpha-\delta_i-\delta_j
\]
for some $1\ls i<j\ls m$.

First suppose $i\ls\len\alpha$. Observe that $\la$ has a $1$-\spam in row $i$, and a $1$-\sprm in row $m$. Since a $2$-Carter partition is necessary \treg, $\alpha_i>\alpha_m$, so the $1$-\spam lies in a longer slope than the $1$-\sprm. So $\noregdown1\la$ by \cref{noregcond}, so $\spn\la$ is inhomogeneous by \cref{atmostev}. 

Now assume that $i>\len\alpha$. In this case we apply a dimension argument. We need to treat the cases $j=m$ and $j<m$ separately.

If $j<m$, then define $\mu=\la+\delta_j-\delta_m$. Then by \cref{dimrowrem,s46} $\la\dblreg=\mu\dblreg$ and $\ddim\la>\ddim\mu$, so $\spn\la$ is inhomogeneous by \cref{samerdoub}.

If $j=m$, define $\mu=\la+\delta_i-\delta_j$. Then \cref{dimrowrem,firstdimlem} imply that $\la\dblreg=\mu\dblreg$ and $\ddim\la>\ddim\mu$, so $\spn\la$ is inhomogeneous.

\subsection*{\nextcase$\la\dn0=\la$, $\la\dn1\in\ird_1$}

In this case we write
\[
\la\dn1=(4m-3\df5,2,1)+4\alpha
\]
where $\alpha$ is a $2$-Carter partition with $\len\alpha\ls m-1$. We can reconstruct $\la$ from this partition by adding nodes of \spr $1$; because $\la\dn0=\la$, we need to add at least one node in every row from $1$ to $m+1$. In each of rows $m,m+1$ we can only add a single $1$-node (so $\la$ ends in $(\dots,3,2)$). So in order to have $\ev\la=2$ we add two nodes in each of rows $1,\dots,m-1$ except one, where we only add one.

In other words,
\[
\la=(4m-1\df7,3,2)+4\alpha-\delta_i
\]
for some $1\ls i\ls m-1$. If $i>\len\alpha$ then $\la$ is \qs, contrary to assumption, so we must have $i\ls\len\alpha$. Now $\la$ has a $1$-\spam in row $i$, and a $1$-\sprm in row $m$ in a shorter slope, and so $\noregdown1\la$ by \cref{noregcond}, and so $\spn\la$ is inhomogeneous by \cref{atmostev}.

\subsection*{\nextcase$\la\dn0=\la$, $\la\dn1\in\ird_2$}

Here
\[
\la\dn1=(4m-3\df5,4,2,1)+4\alpha
\]
where $\alpha$ is a $2$-Carter partition of length at most $m-1$. To reconstruct $\la$ from this with $\la\dn0=\la$, we add:
\begin{itemize}
\item
one or two nodes in each of rows $1,\dots,m-2$;
\item
zero, one or two nodes in row $m-1$ (but we must add at least one node if $\len\alpha=m-1$);
\item
one node in each of rows $m+1,m+2$.
\end{itemize}
In particular, $\la$ ends $(\dots,4,3,2)$. Since $\ev\la=2$, we must therefore add two nodes in each of rows $1,\dots,m-2$, and either zero or two nodes in row $m-1$. If we add two nodes in row $m-1$, then
\[
\la=(4m-1\df7,4,3,2)+4\alpha
\]
which is \qs, contrary to assumption. So instead we must have $\len\alpha\ls m-2$ and
\[
\la=(4m-1\df11,5,4,3,2)+4\alpha.
\]
Now we consider two subcases. Suppose first that $\alpha_1\gs2$. Let $\kappa=\la\up0\up1$. Then we claim that $\noregdown0\kappa$. By \cref{atmostev} this implies that $\spn\kappa$ is inhomogeneous, and hence that $\spn\la$ is inhomogeneous.

We calculate
\begin{align*}
\la\up0&=(4m+1\df13,5,4,3,2,1)+4\alpha,
\\
\intertext{so that}
\kappa=\la\up0\up1&=(4m+3\df15,7,4,3,2,1)+4\alpha.
\end{align*}
$\kappa$ has a $0$-\sprm $(m+3,1)$, lying in slope $2m+4$. In addition, $\kappa$ has a $0$-\spam $(1,4m+4\alpha_1+4)$, lying in slope $2m+2\alpha_1+2$. The assumption that $\alpha_1\gs2$ means that there is a $0$-\spam in a longer slope than the $0$-\sprm, so by \cref{noregcond} $\noregdown0\kappa$, as claimed.

Now we consider the subcase where $\alpha_1\ls1$. Since a $2$-Carter partition is automatically \treg, this means that $\alpha=\varnothing$ or $\alpha=(1)$. Our assumptions then imply that $\len\alpha\ls m-4$, since if $\len\alpha\gs m-3$, then $\la$ is one of the partitions $(5,4,3,2)$, $(11,5,4,3,2)$, $(15,5,4,3,2)$ or $(19,11,5,4,3,2)$. The first of these lies in $\ird_2$, while the others lie in $\scrr$; either way, this contradicts our standing assumptions.

In this situation we consider the partition
\[
\kappa=\la\up0\up1\dn0=(4m+3\df15,7,4,3,2)+4\alpha,
\]
which is \qs. Writing
\[
\kappa=(4m+3\df3)+4(\alpha+(1^{m-2}))\sqcup(4,2),
\]
we see that $\spn\kappa$ is inhomogeneous by \cref{mainqs}, since the partition $\alpha+(1^{m-2})$ is $2$-singular, and so certainly not $2$-Carter. So by \cref{atmostev}, $\spn\la$ is inhomogeneous.

\subsection*{\nextcase$\la\dn0\in\ird_0$, $\la\dn1=\la$}

Here we write
\[
\la\dn0=(4m-1\df3)+4\alpha,
\]
with $\alpha$ a $2$-Carter partition of length at most $m$. To construct $\la$, we add one or two nodes in each of rows $1,\dots,m$, and possibly one node in row $m+1$. So we have two cases.
\begin{enumerate}
\item
$\la=(4m+1\df1)+4\alpha-\delta_i-\delta_j$, where $1\ls i<j\ls m$.

In this case if $i\ls\len\alpha$, then we can use \cref{noregcond} to show that $\noregdown0\la$; this is similar to previous cases. So by \cref{atmostev} $\spn\la$ is inhomogeneous.

So assume instead that $i>\len\alpha$. Now we apply a dimension argument. If we set $\mu=\la+\delta_j-\delta_{m+1}$, then by \cref{dimrowrem,s43} $\la\dblreg=\mu\dblreg$ and $\ddim\la>\ddim\mu$, so that $\spn\la$ is inhomogeneous.

\item
$\la=(4m+1\df5)+4\alpha-\delta_i-\delta_j$, where $1\ls i<j\ls m$.

First we observe that if $i\ls\len\alpha$, then $\noregdown0\la$. This is proved similarly to the case above, and by \cref{atmostev} implies that $\spn\la$ is inhomogeneous.

So assume $i>\len\alpha$. We claim that if $j=m$, then the partition $\kappa=\la\up1$ satisfies $\noregdown0\kappa$, so that $\spn\kappa$, and hence $\spn\la$, are inhomogeneous.

The assumption $j=m$ means that
\begin{align*}
\la&=(4m+1\df4m+9-4i,\ \ 4m+4-4i,\ \ 4m+1-4i\df9,4)+4\alpha
\\
\intertext{and hence}
\kappa=\la\up1&=(4m+3\df4m+11-4i,\ \ 4m+4-4i,\ \ 4m+3-4i\df11,4)+4\alpha.
\end{align*}
$\kappa$ has a $0$-\sprm $(m,4)$ in slope $2m$. Now observe that $m\gs3$, since if $m=2$ then $\la=(8,4)$, contrary to the assumption that $\la\notin\scrr$. So either $i\ls m-2$ or $i\gs2$. In the first case, $\kappa$ has a $0$-\spam $(m-1,12)$; in the other case $\kappa$ has a $0$-\spam $(1,4m+4\alpha_1+4)$. Either way, \cref{noregcond} gives $\noregdown0\kappa$, as claimed.

We are left with the case where $i>\len\alpha$ and $j<m$. Here we apply a dimension argument. Let $\mu=\la+\delta_j-\delta_m$. Then by \cref{dimrowrem,s42} $\la\dblreg=\mu\dblreg$ and $\ddim\la>\ddim\mu$, so $\spn\la$ is inhomogeneous.
\end{enumerate}

\subsection*{\nextcase$\la\dn0\in\ird_1$, $\la\dn1=\la$}

We write
\[
\la\dn0=(4m-1\df3,2)+4\alpha,
\]
where $\alpha$ is $2$-Carter with length at most $m$. To reconstruct $\la$ from $\la\dn0$, we add:
\begin{itemize}
\item
one or two nodes in each of rows $1,\dots,m-1$;
\item
zero, one or two nodes in row $m$ (but we must add at least one node if $\len\alpha=m$);
\item
one node in row $m+2$.
\end{itemize}
Now we consider two cases.
\begin{enumerate}
\item
Suppose we add at least one node in row $m$ when constructing $\la$ from $\la\dn0$. Then
\[
\la=(4m+1\df5,2,1)+4\alpha-\delta_i
\]
for some $i\ls m$. Furthermore, we have $i\ls\len\alpha$, since if $i>\len\alpha$ then $\la$ is \qs, contrary to assumption. Now we claim that in most cases $\noregdown0\la$, so that $\la$ is inhomogeneous  by \cref{atmostev}.

If $i<m$ then $\la$ has a $0$-\sprm $(m,5+4\alpha_m)$ and a $0$-\spam $(i,4m+5-4i+4\alpha_i)$. The fact that $\alpha_i>\alpha_m$ then gives $\noregdown0\la$, by \cref{noregcond}.

If $i=m$ and $\alpha_m\gs2$, then $\la$ has a $0$-\sprm $(m+2,1)$ and a $0$-\spam $(m,5+4\alpha_m)$, and again \cref{noregcond} gives $\noregdown0\la$.

So we are done unless $i=m$ and $\alpha_m=1$. In this case, we let $\kappa=\la\up1$, and we claim that $\noregdown0\kappa$, so that $\spn\kappa$ (and hence $\spn\la$) are inhomogeneous.

We calculate
\[
\kappa=(4m+3\df11,4,3,2)+4\alpha.
\]
This has a $0$-\sprm $(m,8)$ and a $0$-\spam $(m+3,1)$, and by \cref{noregcond} we have $\noregdown0\kappa$, as claimed.
\item
Now suppose that in reconstructing $\la$ from $\la\dn0$ we do not add a node in row $m$. Then
\[
\la=(4m+1\df9,3,2,1)+4\alpha-\delta_i,
\]
with $\len\alpha\ls m-1$ and $i\ls m-1$. Again, we claim that in most cases $\noregdown0\la$.

If $\alpha_i\gs2$, then $\la$ has a $0$-\sprm $(m+2,1)$ and a $0$-\spam $(i,4m+5-4i+4\alpha_i)$, and \cref{noregcond} gives $\noregdown0\la$. If $\alpha_i=1$ and $i\ls m-2$, then $i=\len\alpha$, so $\la$ has a $0$-\sprm $(m-1,9)$ and a $0$-\spam $(i,4m+9-4i)$, and again \cref{noregcond} applies.

We split the remaining possibilities into two cases.
\begin{enumerate}
\item
Suppose $i=m-1$ and $\alpha_i\ls1$. Then
\begin{align*}
\la&=(4m+1\df13,8,3,2,1)+4\alpha.
\\
\intertext{Now set $\kappa=\la\up1$; then we claim that $\noregdown0\kappa$. We have}
\kappa&=(4m+3\df15,8,3,2,1)+4\alpha.
\end{align*}
$m$ must be greater than or equal to $3$ here, since if $m=2$ then $\la=(8,3,2,1)$ or $(12,3,2,1)$, contrary to our standing assumption that $\la\notin\scrr$. But this means that $\kappa$ has a $0$-\spam $(1,4m+4\alpha_1+4)$ and a $0$-\sprm $(m-1,8)$, so \cref{noregcond} gives $\noregdown0\kappa$.

\item
Now suppose $\len\alpha<i<m-1$. Consider the partition
\[
\kappa=\la\up1\dn0\dn1\dn0=(4m-1\df7,2)+4\alpha.
\]
Observe that $\kappa$ is \qs, and write
\[
\kappa=(4m-5\df3)+4(\alpha+(1^{m-1}))\sqcup2.
\]
Because $\len\alpha\ls m-3$, the partition $\alpha+(1^{m-1})$ is $2$-singular and so certainly not $2$-Carter, and so by \cref{mainqs} $\spn\kappa$ is inhomogeneous. So by four applications of \cref{atmostev}, $\spn\la$ is inhomogeneous.

\end{enumerate}
\end{enumerate}

\subsection*{\nextcase$\la\dn0\in\ird_2$, $\la\dn1=\la$}

We write
\[
\la\dn0=(4m-1\df7,4,3,2)+4\alpha
\]
where $\alpha$ is a $2$-Carter partition of length at most $m-1$. To reconstruct $\la$ we add nodes of \spr $0$. To ensure that $\la$ has no $1$-\sprms, we need to add
\begin{itemize}
\item
either one or two nodes in each of rows $1,\dots,m-1$,
\item
at most one node in row $m$, and
\item
one node in row $m+3$.
\end{itemize}
(Note that we cannot add a node in row $m+1$ because the condition that $\ev\la=2$ would force $\la=(4m+1\df9,5,4,2,1)+4\alpha$, which is \qs, contrary to assumption.)

The assumption that $\ev\la=2$ gives
\[
\la=(4m+1\df5,3,2,1)+4\alpha-\delta_i
\]
for some $1\ls i\ls m$. We must have $m\gs2$, since if $m=1$ then $\la=(4,3,2,1)$, contrary to assumption.

We claim that in most cases $\noregdown0\la$. First observe that $\la$ has a $0$-\spam $(m+1,4)$. If either $\len\alpha\ls m-2$ or $i\ls m-1$, then $\la$ also has either $(m-1,9)$ or $(m,5)$ as a $0$-\sprm, in which case \cref{noregcond} gives $\noregdown0\la$.

We are left with the case where $\len\alpha=m-1$ and $i=m$.  Now we must have $\alpha_1\gs2$, since otherwise the fact that $\alpha$ is \treg would give $\la=(13,4,3,2,1)$, contrary to our standing assumption that $\la\notin\scrr$. We define

\[
\kappa=\la\up1=(4m+3,\dots,11,4,3,2,1)+4\alpha.
\]
Observe that $\kappa$ has a $0$-\sprm $(m+3,1)$ and a $0$-\spam $(1,4m+4+4\alpha_1)$. So $\noregdown0\kappa$ by \cref{noregcond}, and hence $\spn\kappa$ and $\spn\la$ are inhomogeneous.

\subsection*{\nextcase$\la\dn0\in\ird$, $\la\dn1\in\ird$}

In this case $\la_r=\max\{(\la\dn0)_r,(\la\dn1)_r\}$ for every $r$, because if $i$ denotes the \spr of the node at the end of row $r$ of $\la$, then $(\la\dn{1-i})_r=\la_r\gs(\la\dn i)_r$. The only positive even parts that $\la\dn0$ or $\la\dn1$ can have are $2$ and $4$, so the positive even parts of $\la$ are $2$ and $4$. So either $\la\dn0$ or $\la\dn1$ contains $4$ as a part, and hence lies in $\ird_2$. We consider the possibilities.

\begin{description}
\item[$\la\dn0\in \ird_2$, $\la\dn1\in\ird_0\cup\ird_1${\rm:}]
Here
\[
\la\dn0=(4m+3\df7,4,3,2)+4\alpha,\qquad\la\dn1=(4n+1\df5,[2,]1)+4\beta
\]
for some $2$-Carter partitions $\alpha,\beta$ with $\len\alpha\ls m$ and $\len\beta\ls n+1$; the notation $[2,]$ indicates that the part $2$ may or may not be present in $(\la\dn1)$. In order for $\la$ to contain $4$ as a part, we need $(\la\dn1)_{m+1}\ls4$, so that $m\gs n$. But then $\la_{m+3}=2$ while $(\la\dn1)_{m+3}=0$, which is impossible.

\item[$\la\dn0\in\ird_0$, $\la\dn1\in\ird_2${\rm:}]
Here
\[
\la\dn0=(4m+3\df7,3)+4\alpha,\qquad\la\dn1=(4n+1\df5,4,2,1)+4\beta.
\]
Now in order to have $4$ as a part of $\la$ we need $(\la\dn0)_{n+1}\ls4$, giving $m\ls n$. But now $\la_{n+2}=2$ while $(\la\dn0)_{n+2}=0$, which is impossible.

\item[$\la\dn0\in\ird_1$, $\la\dn1\in\ird_2${\rm:}]
Now
\[
\la\dn0=(4m+3\df7,3,2)+4\alpha,\qquad\la\dn1=(4n+1\df5,4,2,1)+4\beta.
\]
Now the only way to have $4$ as a part of $\la$ is if $m=n$, in which case the last three parts of $\la$ are $(\la_{m+1},\la_{m+2},\la_{m+3})=(4,2,1)$. Since $\la\neq\la\dn1$, there must be some $r\ls m$ such that $\la_r=(\la\dn0)_r=4m+7-4r+4\alpha_r$. If $\alpha_r=0$ this means that $\la$ has a $1$-\sprm $(r,\la_r)$ and a $1$-\spam $(m+3,2)$ in a longer slope, so that $\noregdown1\la$. If $\alpha_r>0$, then $\la$ has a $0$-\spam $(r,\la_r+1)$ and a $0$-\sprm $(m+1,4)$ in a shorter slope, so $\noregdown0\la$. Either way, $\spn\la$ is inhomogeneous by \cref{atmostev}.

\item[$\la\dn0,\la\dn1\in\ird_2${\rm:}]
Now
\[
\la\dn0=(4m+3\df7,4,3,2)+4\alpha,\qquad\la\dn1=(4n+1\df5,4,2,1)+4\beta.
\]
To get $4$ as a part of $\la$, we must have either $m=n-1$ or $m=n$. If $m=n-1$, then the last three non-zero parts of $\la$ are $(4,2,1)$, and we can proceed as in the previous case. If $m=n$, then the last three parts of $\la$ are $(\la_{m+1},\la_{m+2},\la_{m+3})=(4,3,2)$. Since $\la\neq\la\dn0$, there must be some $r\ls m$ such that $\la_r=(\la\dn1)_r=4m+5-4r+4\beta_r$. If $\beta_r\ls1$, then $\la$ has a $0$-\sprm $(r,\la_r)$ and a $0$-\spam $(m+4,1)$ in a longer slope, so $\noregdown0\la$. If $\beta_r\gs2$, then $\la$ has a $1$-\spam $(r,\la_r+1)$ and a $1$-\sprm $(m+2,3)$ is a shorter slope, so $\noregdown1\la$. Either way, $\spn\la$ is inhomogeneous.
\end{description}

\begin{pf}[Proof of \cref{mainthm} (``only if'' part)]
We proceed by induction on $\card\la$. Suppose $\la\in\cald$ with $\ev\la=2$ and $\la\notin\scri_2$, and that the theorem is known to be true for all \trps smaller than $\la$. If $\la$ is separated, then \cref{mainqs} gives the result. If $\la$ is one of the seven partitions in $\scrr$, then \cref{scrrred} gives the result. If for $i=0$ or $1$ the partition $\la\dn i$ does not lie in $\scri\cup\{\la\}$, then \cref{atmostev} together with the inductive hypothesis gives the result. Otherwise, $\la$ satisfies the assumptions set out following \cref{atmostev}, and so is dealt with by Cases 1--\arabic{case} above.
\end{pf}

\section{Index of notation}

For the reader's convenience we conclude with an index of the notation we use in this paper. We provide references to the relevant subsections.

\newlength\colwi
\newlength\colwii
\newlength\colwiii
\setlength\colwi{2.5cm}
\setlength\colwiii{1cm}
\setlength\colwii\textwidth
\addtolength\colwii{-\colwi}
\addtolength\colwii{-\colwiii}
\addtolength\colwii{-1em}
\subsubsection*{Partitions}
\vspace{-\topsep}
\begin{longtable}{@{}p{\colwi}p{\colwii}p{\colwiii}@{}}
$\calp$&the set of all partitions&\ref{compnpartnsec}\\
$\calp(n)$&the set of all partitions of $n$&\ref{compnpartnsec}\\
$\cald$&the set of all \trps&\ref{compnpartnsec}\\
$\cald(n)$&the set of all \trps of $n$&\ref{compnpartnsec}\\
$\varnothing$&the partition of $0$&\ref{compnpartnsec}\\
$\len\la$&the length of a partition $\la$&\ref{compnpartnsec}\\
$\ev\la$&the number of positive even parts of a partition $\la$&\ref{snsec}\\
$\la'$&the partition conjugate to $\la$&\ref{compnpartnsec}\\
$\dom$&the dominance order on $\calp(n)$&\ref{compnpartnsec}\\
$a\la$&the partition $(a\la_1,a\la_2,\dots)$&\ref{compnpartnsec}\\
$\la+\mu$&the partition $(\la_1+\mu_1,\la_2+\mu_2,\dots)$&\ref{compnpartnsec}\\
$\la\sqcup\mu$&the partition obtained by arranging all the parts of $\la$ and $\mu$ together in decreasing order&\ref{compnpartnsec}\\
$a\df b$&the arithmetic progression $a,a-4,\dots,b$&\ref{compnpartnsec}\\
$\epstar\mu$&the $2$-sign of $\mu'$ if $\mu$ has $2$-core $\varnothing$&\ref{corequotsec}\\
$\la\reg$&the regularisation of $\la\in\calp(n)$&\ref{regdoubsec}\\
$\la\dbl$&the double of $\la\in\cald(n)$&\ref{regdoubsec}\\
$\la\dblreg$&$(\la\dbl)\reg$&\ref{regdoubsec}\\
$\lrc\gamma\alpha\beta$&the Littlewood--Richardson coefficient corresponding to $\alpha,\beta,\gamma\in\calp$&\ref{lrsec}\\
$\lrcb\gamma\alpha$&$\sum_{a\gs0}\lrc\gamma\alpha{(1^a)}$&\ref{lrsec}\\
\end{longtable}

\subsubsection*{Symmetric functions}
\vspace{-\topsep}
\begin{longtable}{@{}p{\colwi}p{\colwii}p{\colwiii}@{}}
$\La$&the algebra of symmetric functions&\ref{algsymfnsec}\\
$\Delta$&coproduct on $\La$&\ref{algsymfnsec}\\
$m_\la$&the monomial symmetric function corresponding to $\la\in\calp$&\ref{algsymfnsec}\\
$\sstd\la$&the set of semistandard $\la$-tableaux with entries in $X$&\ref{algsymfnsec}\\
$s_\la$&the Schur function corresponding to $\la\in\calp$&\ref{algsymfnsec}\\
$\shtd\la$&the set of semistandard shifted $\la$-tableaux&\ref{algsymfnsec}\\
$P_\la$&the Schur P-function corresponding to $\la\in\cald$&\ref{algsymfnsec}\\
$\ip{\ \ }{\ \ }$&standard inner product on $\sym$&\ref{algsymfnsec}\\
$\omega$&involution on $\La$ defined by $s_\la\mapsto s_{\la'}$&\ref{algsymfnsec}\\
$\ddd_\nu$&reduction operator on $\La$&\ref{redopsec}\\
$\strp\la\mu$&integer depending on two \trps $\la,\mu$&\ref{redopsec}\\
\end{longtable}

\subsubsection*{Groups, algebras and representations}
\vspace{-\topsep}
\begin{longtable}{@{}p{\colwi}p{\colwii}p{\colwiii}@{}}
$\sss n$&the symmetric group of degree $n$&\ref{snsec}\\
$\aaa n$&the alternating group of degree $n$&\ref{ansec}\\
$\tsss n$&a double cover of $\sss n$&\ref{snsec}\\
$\taaa n$&a double cover of $\aaa n$&\ref{ansec}\\
$\schur n$&the Schur algebra of degree $n$ over $\bbf$&\ref{schursec}\\
$\qschur n$&the $q$-Schur algebra of degree $n$ over $\bbc$, with $q=-1$&\ref{schursec}\\
$\spe\la$&the Specht module for $\sss n$ corresponding to $\la\in\calp(n)$&\ref{snsec}\\
$\jms\la$&the James module for $\sss n$ corresponding to $\la\in\cald(n)$&\ref{snsec}\\
$\weyl\la$&the Weyl module for $\schur n$ or $\qschur n$ corresponding to $\la\in\calp(n)$&\ref{schursec}\\
$\schir\la$&the irreducible module for $\schur n$ or $\qschur n$ corresponding to $\la\in\calp(n)$&\ref{schursec}\\
$D_{\la\mu}$&the decomposition number $\cm{\weyl\la}{\schir\mu}$ for $\schur n$ (equals the decomposition number $\cm{\spe\la}{\jms\mu}$ for $\bbf\sss n$ if $\mu\in\cald(n)$)&\ref{decnosec}\\
$\odc_{\la\mu}$&the decomposition number $\cm{\weyl\la}{\schir\mu}$ for $\qschur n$&\ref{schursec}\\
$A_{\la\mu}$&the $(\la,\mu)$ entry of the adjustment matrix for $\schur n$&\ref{schursec}\\
$\ord\la$&the character of $\spe\la$ over $\bbc$&\ref{snsec}\\
$\sid\la$&the Brauer character of $\jms\la$&\ref{snsec}\\
$\spn\la$&a class function labelled by $\la\in\cald(n)$&\ref{snsec}\\
$\chm{\ \ }{\ \ }$&the standard inner product on characters&\ref{snsec}\\
$\decs\la\mu$&the modified decomposition number $\cm{\mspn\la}{\sid\mu}$&\ref{decnosec}\\
$\prj\mu$&the character of the projective cover of $\jms\mu$&\ref{rouqdecsec}\\
$\omega^\mu$&a virtual projective character in a Rouquier block&\ref{rouqdecsec}\\
$\dim\la$&the degree of the character $\spn\la$&\ref{degsec}\\
$\ddim\la$&the divided degree of $\spn\la$&\ref{degsec}\\
$\scrr$&a set of seven exceptional partitions&\ref{mainproofsec}\\
$\ird$&the set of \trps with at most two even parts labelling homogeneous spin characters&\ref{mainproofsec}\\
$\ird_k$&the set of partitions in $\ird$ with $k$ even parts&\ref{mainproofsec}\\
\end{longtable}

\needspace{5em}
\subsubsection*{Branching rules}
\vspace{-\topsep}
\begin{longtable}{@{}p{\colwi}p{\colwii}p{\colwiii}@{}}
$\chi\res_{\tsss{n-1}}$&the restriction of $\chi$ to $\tsss{n-1}$&\ref{brnchsec}\\
$\chi\ind^{\tsss{n+1}}$&the character obtained by inducing $\chi$ to $\tsss{n+1}$&\ref{brnchsec}\\
$\ee i$&Robinson's $i$-restriction functor&\ref{brnchsec}\\
$\ff i$&Robinson's $i$-induction functor&\ref{brnchsec}\\
$\eed ir$&$\ee i^r/r!$&\ref{brnchsec}\\
$\ffd ir$&$\ff i^r/r!$&\ref{brnchsec}\\
$\epsilon_i\chi$&$\max\lset{r\gs0}{\ee i^r\chi\neq0}$&\ref{brnchsec}\\
$\varphi_i\chi$&$\max\lset{r\gs0}{\ff i^r\chi\neq0}$&\ref{brnchsec}\\
$\emx i\chi$&$\eed i{\epsilon_i\chi}\chi$&\ref{brnchsec}\\
$\fmx i\chi$&$\ffd i{\varphi_i\chi}\chi$&\ref{brnchsec}\\
$\la\dn i$&the partition obtained by removing all the $i$-\sprms of $\la$&\ref{partnsec}\\
$\la\up i$&the partition obtained by adding all the $i$-\spams of $\la$&\ref{partnsec}\\
$\noregdown i\la$&$\epsilon_i\sid{\la\dblreg}<\epsilon_i\spn\la$&\ref{brnchsec}\\
\end{longtable}

\end{document}